\documentclass[10pt, a4paper, reqno]{amsart}

\usepackage{amsthm,amsmath, amsfonts, amssymb, amscd,enumerate,enumitem, amscd, esint, xcolor, url, mathrsfs,}
\usepackage[colorlinks]{hyperref}
\usepackage[bbgreekl]{mathbbol}
\usepackage[utf8]{inputenc}
\usepackage{mathtools}

\DeclareSymbolFontAlphabet{\amsbb}{bbold}
\DeclareSymbolFontAlphabet{\mathbb}{AMSb}%

\newtheorem{thm}{Theorem}[section]

\newtheorem{prop}[thm]{Proposition}
\theoremstyle{definition}

\theoremstyle{remark}
\newtheorem{rem}[thm]{Remark}

\DeclareMathOperator{\supp}{supp}

\DeclareMathOperator*{\essinf}{ess\ inf}

\newcommand{\zinfa}{(0,\infty)\times (0,\infty)\times (-1,1)}
\newcommand{\loc}{\textup{loc}}
\newcommand{\glob}{\textup{glob}}

\numberwithin{equation}{section}
\usepackage[defaultlines=4, all]{nowidow}
\setnoclub[2]
\allowdisplaybreaks[3]

\setlist[enumerate,1]{label=(\alph*)}

\begin{document}
\raggedbottom %Delete ‘Underfull \vbox (badness 10000) has occurred while \output is active'

    \title[BLO spaces associated with Laguerre polynomials expansions]
    {BLO spaces associated with Laguerre polynomial expansions}

    \author[J. J. Betancor]{Jorge J. Betancor}
    \address{Jorge J. Betancor\newline
        Departamento de Análisis Matemático, Universidad de La Laguna,\newline
        Campus de Anchieta, Avda. Astrofísico Sánchez, s/n,\newline
        38721 La Laguna (Sta. Cruz de Tenerife), Spain}
    \email{jbetanco@ull.es}

    \author[E. Dalmasso]{Estefanía Dalmasso}
    \address{Estefanía Dalmasso, Pablo Quijano\newline
        Instituto de Matemática Aplicada del Litoral, UNL, CONICET, FIQ.\newline Colectora Ruta Nac. Nº 168, Paraje El Pozo,\newline S3007ABA, Santa Fe, Argentina}
    \email{edalmasso@santafe-conicet.gov.ar, pquijano@santafe-conicet.gov.ar}

    \author[P. Quijano]{Pablo Quijano}

    \thanks{The first author is partially supported by grant PID2019-106093GB-I00 from the Spanish Government. The second author is partially supported by grants PICT-2019-2019-00389 (ANPCyT), PIP-1220200101916O (CONICET) and CAI+D 2019-015 (UNL)}
    \date{\today}


    %%% ----------------------------------------------------------------------
    \begin{abstract}
     In this paper we introduce spaces of $\textup{BLO}$-type related to Laguerre polynomial expansions. We consider the probability measure on $(0,\infty)$ defined by $d\gamma_\alpha(x)=\frac{2}{\Gamma(\alpha+1)}e^{-x^2}x^{2\alpha+1}dx$ with $\alpha>-\frac12$. For every $a>0$, the space $\textup{BLO}_a((0,\infty),\gamma_\alpha)$ consists of all those measurable functions defined on $(0,\infty)$ having bounded lower oscillation with respect to $\gamma_\alpha$ over an admissible family $\mathcal{B}_a$ of intervals in $(0,\infty)$. The space $\textup{BLO}_a((0,\infty),\gamma_\alpha)$ is a subspace of the space $\textup{BMO}_a((0,\infty),\gamma_\alpha)$ of bounded mean oscillation functions with respect to $\gamma_\alpha$ and $\mathcal{B}_a$. The natural $a$-local centered maximal function defined by $\gamma_\alpha$ is bounded from $\textup{BMO}_a((0,\infty),\gamma_\alpha)$ into $\textup{BLO}_a((0,\infty),\gamma_\alpha)$. We prove that the maximal operator, the $\rho$-variation and the oscillation operators associated with local truncations of the Riesz transforms in the Laguerre setting are bounded from $L^\infty((0,\infty),\gamma_\alpha)$ into $\textup{BLO}_a((0,\infty),\gamma_\alpha)$. Also, we obtain a similar result for the maximal operator of local truncations for spectral Laplace transform type multipliers.
    \end{abstract}
    %%% ----------------------------------------------------------------------
    \maketitle
    %%% ----------------------------------------------------------------------

\section{Introduction}

We consider, for every $\alpha>-\frac12$, the probability measure defined on $(0,\infty)$ by $d\gamma_\alpha(x)=\frac{2}{\Gamma(\alpha+1)}e^{-x^2}x^{2\alpha+1}dx$. This measure has not the doubling property with respect to the usual metric defined by the absolute value $|\cdot|$ on $(0,\infty)$. Then, the triple $((0,\infty),|\cdot|,\gamma_\alpha)$ is not homogeneous in the sense of Coifman and Weiss (\cite{CW}). Harmonic analysis in the spaces of homogeneous type can be developed following the model of Euclidean spaces $(\mathbb{R}^n,\|\cdot\|,\lambda)$ where $\|\cdot\|$ denotes a norm and $\lambda$ is the Lebesgue measure on $\mathbb{R}^n$. When the measure is not doubling the situation is very different and it is necessary to introduce new ideas (see, for instance, \cite{CMY}, \cite{FLYY}, \cite{FYY}, \cite{Hy}, \cite{HyYY}, \cite{NTV}, \cite{To1}, \cite{To2} and \cite{To3}).

Tolsa (\cite{To1}) defined $\textup{BMO}$-type spaces, that he named $\textup{RBMO}$-spaces, on $(\mathbb{R}^n,\mu)$ when $\mu$ is a Radon measure on $\mathbb{R}^n$, which is not necessarily doubling, satisfying that $\mu(B(x,r))\le Cr^k$, $x\in \mathbb{R}^k$ and $r>0$, for some $k\in \{1,\dots,n\}$ and $C>0$. He also proved that $\textup{RBMO}(\mathbb{R}^n,\mu)$ has many of the properties of the classical space $\textup{BMO}(\mathbb{R}^n)$ of John and Nirenberg. In particular, the integral operators defined by standard Calderón-Zygmund kernels are bounded from $L^\infty(\mathbb{R}^n,\mu)$ into $\textup{RBMO}(\mathbb{R}^n,\mu)$.

It is clear that, for every $0<r\le x$, $\gamma_\alpha((x-r,x+r))\le Cr$. Then, following Tolsa's ideas we can define the space $\textup{RBMO}((0,\infty),\gamma_\alpha)$ by replacing $\mathbb{R}^n$ by $(0,\infty)$. However, $\textup{RBMO}((0,\infty),\gamma_\alpha)$ is not suitable to study harmonic analysis operators associated with Laguerre polynomial expansions because these operators are not defined by standard Calderón-Zygmund kernels (\cite{FSS1}, \cite{SaSpec} and \cite{SaFunct}). Motivated by the results in \cite{MM} in the Gaussian setting, the authors and R. Scotto (\cite{BDQS1}) defined a local $\textup{BMO}$-type space related to the measure $\gamma_\alpha$ as follows.

We consider the function $m(x)=\min\{1,1/x\}$, $x\in (0,\infty)$. Given $a>0$, we say that an interval $(x-r,x+r)$, with $0<r\le x$, is \textit{$a$-admissible}, or is in the class $\mathcal{B}_a$, when $r\le am(x)$. The measure $\gamma_\alpha$ has the doubling property on $\mathcal{B}_a$, that is, there exists $C>0$ such that, for every $0<r\le x$ being $r\le am(x)$, we have that
\[
\gamma_\alpha(I(x,2r))\le C\gamma_\alpha((x-r,x+r)),
\]
where $I(x,r):=(x-r,x+r)\cap (0,\infty)$ for $x,r>0$.

A function $f\in L^1((0,\infty),\gamma_\alpha)$ is said to be in $\textup{BMO}_a((0,\infty),\gamma_\alpha)$ when
\[
\|f\|_{*,\alpha,a}:=\sup_{I\in \mathcal{B}_a}\frac{1}
{\gamma_\alpha(I)}\int_I|f(y)-f_I|d\gamma_\alpha(y)<\infty,\]
where $f_I=\fint_If(y)d\gamma_\alpha(y)$, for every $I\in \mathcal{B}_a$. For every $f\in \textup{BMO}_a((0,\infty),\gamma_\alpha)$, we define
\[\|f\|_{\textup{BMO}_a((0,\infty),\gamma_\alpha)}:=\|f\|_{L^1((0,\infty),\gamma_\alpha)}+\|f\|_{*,\alpha,a}.
\]
The space $\textup{BMO}_a((0,\infty),\gamma_\alpha)$ actually does not depend on $a>0$. Then, in the sequel we will write $\textup{BMO}((0,\infty),\gamma_\alpha)$ and $\|\cdot\|_{*,\alpha}$ instead of $\textup{BMO}_a((0,\infty),\gamma_\alpha)$ and $\|\cdot\|_{*,\alpha,a}$, respectively. This space can be identified with the dual space of the Hardy space $H^1((0,\infty),\gamma_\alpha)$ studied in \cite{BDQS1} (see \cite[Theorem~1.1]{BDQS1}).

The space $\textup{BLO}(\mathbb{R}^n)$ of functions of bounded lower oscillation on $\mathbb{R}^n$ was introduced by Coifman and Rochberg (\cite{CR}). Later, Bennett (\cite{Be}) obtained a characterization of the functions in $\textup{BLO}(\mathbb{R}^n)$ by using the natural Hardy-Littlewood maximal operators, and Leckband (\cite{Lec}) proved that certain maximal operators associated with singular integrals are bounded from $L^{p}(\mathbb{R}^n)\cap L^\infty(\mathbb{R}^n)$ into $\textup{BLO}(\mathbb{R}^n)$ for certain $1\le p<\infty$.

Based on Tolsa's ideas, Jiang (\cite{Ji}) introduced $\textup{BLO}$-type spaces in $(\mathbb{R}^n,\mu)$ where~$\mu$ is a positive non-doubling Radon measure with polynomial growth. $\textup{BLO}$-spaces in the Gaussian setting were defined by Liu and Yang (\cite{LiuY}). In \cite{HyYY}, Littlewood-Paley functions in non-doubling settings on $\textup{RBLO}$ spaces were studied. Other results concerning $\textup{RBLO}$ spaces can be encountered in \cite{LLW} and \cite{LinY}. As in happens with $\textup{RBMO}$-spaces, $\textup{RBLO}$-spaces for $\gamma_\alpha$ do not work in a correct way in connection with harmonic analysis operators associated to Laguerre polynomial expansions. 

In this paper we introduce $\textup{BLO}$-spaces associated with the measure $\gamma_\alpha$ on $(0,\infty)$ by using admissible intervals. 

Let $a>0$. We say that a function $f\in L^1((0,\infty),\gamma_\alpha)$ is in $\textup{BLO}_a((0,\infty),\gamma_\alpha)$ when
\[\sup_{I\in \mathcal{B}_a}\frac{1}{\gamma_\alpha(I)}\int_I \left(f(y)-\essinf_{z\in I}f(z)\right)d\gamma_\alpha(y)<\infty.
\]
For every $f\in \textup{BLO}_a((0,\infty),\gamma_\alpha)$ we define
\[\|f\|_{\textup{BLO}_a((0,\infty),\gamma_\alpha)}:=\|f\|_{L^1((0,\infty),\gamma_\alpha)}+\sup_{I\in \mathcal{B}_a}\frac{1}{\gamma_\alpha(I)}\int_I \left(f(y)-\essinf_{z\in I}f(z)\right)d\gamma_\alpha(y).\]
It is not hard to see that
\[
L^\infty((0,\infty),\gamma_\alpha)\subset \textup{BLO}_a((0,\infty),\gamma_\alpha)\subset \textup{BMO}_a((0,\infty),\gamma_\alpha).
\]
The main properties of the space $\textup{BLO}_a((0,\infty),\gamma_\alpha)$ will be established in Section~\ref{sec: BLO}.

Our objective is to prove that maximal, variation and oscillation operators defined by singular integrals in the Laguerre settings are bounded from $L^\infty((0,\infty),\gamma_\alpha)$ to $\textup{BLO}_a((0,\infty),\gamma_\alpha)$.

We now define the operators we are going to consider. Let $\alpha>-\frac12$. The Laguerre polynomial $L_k^\alpha$ of order $\alpha$ and degree $k\in \mathbb{N}$ (see \cite{Leb}) is
\[
L_k^\alpha(x)=\sqrt{\frac{\Gamma(\alpha+1)}{\Gamma(\alpha+k+1)k!}}e^x x^{-\alpha}\frac{d^k}{dx^k}(e^{-x}x^{\alpha+k}),\quad x\in (0,\infty).
\]
The Laguerre differential operator $\widetilde{\Delta_\alpha}$ is given by
\[
\widetilde{\Delta_\alpha}:=\frac{1}{2}\frac{d^2}{dx^2}+\left(\frac{2\alpha+1}{2x}-x\right)\frac{d}{dx}+\alpha+1, \quad f\in C^2(0,\infty).
\]
We define, for every $k\in \mathbb{N}$, $\mathcal{L}_k^\alpha(x):=L_k^\alpha(x^2)$, $x\in (0,\infty)$. Then, the sequence $\{\mathcal{L}_k^\alpha\}_{k\in \mathbb{N}}$ is an orthonormal basis on $L^2((0,\infty),\gamma_\alpha)$. For every $k\in \mathbb{N}$, $\mathcal{L}_k^\alpha$ is an eigenfunction for $\widetilde{\Delta_\alpha}$ associated with the eigenvalue $\lambda_k^\alpha=2k+\alpha+1$.

For every $f\in L^1((0,\infty),\gamma_\alpha)$, we define
\[
c_k^\alpha(f):=\int_0^\infty f(y)\mathcal{L}_k^\alpha(x)d\gamma_\alpha(x),\quad k\in \mathbb{N}.
\]
We consider the operator $\Delta_\alpha$ given by
\[
\Delta_\alpha f=\sum_{k=0}^\infty \lambda_k^\alpha c_k^\alpha(f)\mathcal{L}_k^\alpha,\quad f\in D(\Delta_\alpha),
\]
being
\[
D(\Delta_\alpha)=\left\{f\in L^2((0,\infty),\gamma_\alpha):\,\sum_{k=0}^\infty (\lambda_k^\alpha|c_k^\alpha(f)|)^2<\infty\right\}.
\]
The space $C_c^\infty(0,\infty)$ of all the smooth functions with compact support in $(0,\infty)$ is contained in $D(\Delta_\alpha)$ and $\Delta_\alpha f=\widetilde{\Delta_\alpha}f$, for any $f\in C_c^\infty(0,\infty)$. The operator $\Delta_\alpha$ is self-adjoint and positive in $L^2((0,\infty),\gamma_\alpha)$. Furthermore, the operator $-\Delta_\alpha$ generates a $C_0$-semigroup of operators $\{W_t^{\alpha
}\}_{t>0}$, where, for every $t>0$,
\[
W_t^\alpha(f)=\sum_{k=0}
^\infty e^{-\lambda_k^\alpha t}c_k^\alpha(f)\mathcal{L}^\alpha_k,\quad f\in L^2((0,\infty),\gamma_\alpha).
\]
According to \cite[(4.17.6)]{Leb} we have that, for every $x,y,t\in (0,\infty)$,
\begin{align}\label{1.1}
\nonumber\sum_{k=0}^\infty e^{-kt}&\mathcal{L}_k^\alpha(x)\mathcal{L}_k^\alpha(y)\\
&=\frac{\Gamma(\alpha+1)}{1-e^{-t}}(e^{-t/2}xy)^{-\alpha}I_\alpha\left(\frac{2e^{-t/2}xy}{1-e^{-t}}\right)\exp\left(-\frac{e^{-t}(x^2+y^2)}{1-e^{-t}}\right),
\end{align}
being $I_\alpha$ the modified Bessel function of the first kind and order $\alpha$.

By using \eqref{1.1} we can write, for every $f\in L^2((0,\infty),\gamma_\alpha)$ and $t>0$,
\begin{equation}\label{1.2}
W_t^\alpha(f)(x)=\int_0^\infty W_t^\alpha(x,y)f(y)d\gamma_\alpha(y),
\quad x\in (0,\infty),
\end{equation}
where
\[W_t^\alpha(x,y)=\frac{\Gamma(\alpha+1)e^{-t(\alpha+1)}}{1-e^{-2t}}(e^{-t}xy)^{-\alpha}I_\alpha\left(\frac{2e^{-t}xy}{1-e^{-2t}}\right)\exp\left(-\frac{e^{-2t}(x^2+y^2)}{1-e^{-2t}}\right), \]
for $x,y,t\in (0,\infty).$

The integral in \eqref{1.2} is absolutely convergent for every $f\in L^p((0,\infty),\gamma_\alpha)$, ${1\le p<\infty}$, and for every $t,x\in (0,\infty)$. By defining $W_t^\alpha(f)$ by \eqref{1.2}, for every $f\in L^p((0,\infty),\gamma_\alpha)$ and $t>0$, the family $\{W_t^\alpha\}_{t>0}$ is a $C_0$-semigroup in $L^p((0,\infty),\gamma_\alpha)$, for every $1\le p<\infty$. Thus $\{W_t^\alpha\}_{t>0}$ is a symmetric diffusion semigroup in the sense of Stein (\cite{StLP}).

The study of harmonic analysis in Laguerre settings was begun by Muckenhoupt (\cite{Mu2}) who proved that the maximal operator $W_*^\alpha$ defined by
\[W_*^\alpha(f)=\sup_{t>0}|W_t^\alpha(f)|\]
is bounded from $L^1((0,\infty),\gamma_\alpha)$ into $L^{1,\infty}((0,\infty,\gamma_\alpha)$.  This property was generalized by Dinger (\cite{Di}) to higher dimensions.

We define the Riesz transform $R^\alpha$ associated with the Laguerre operator $\Delta_\alpha$ by
\[R^\alpha(f)=\sum_{k=1}^\infty\frac{1}{\sqrt{\lambda_k^\alpha}}c_k^\alpha(f)\frac{d}{dx}\mathcal{L}_k^\alpha,\quad f\in L^2((0,\infty),\gamma_\alpha).\]
Thus $R^\alpha$ defines a bounded operator on $L^2((0,\infty),\gamma_\alpha)$ (see \cite{NS06}). Furthermore, $R^\alpha$ can be extended from $L^2((0,\infty),\gamma_\alpha)\cap L^p((0,\infty),\gamma_\alpha)$ as a bounded operator on $L^p((0,\infty),\gamma_\alpha)$, for every $1<p<\infty$, and from $L^1((0,\infty),\gamma_\alpha)$ into $L^{1,\infty}((0,\infty),\gamma_\alpha)$ (\cite{SaWeak}). The authors and R. Scotto (\cite{BDQS2}) extended the above results by considering variable exponents $L^{p(\cdot)}$-spaces. Also, in \cite{BDQS1}, endpoint estimates for Riesz transform $R^\alpha$ were established proving that $R^\alpha$ defines a bounded operator from $H^1((0,\infty),\gamma_\alpha)$ into $L^1((0,\infty),\gamma_\alpha)$ and from $L^\infty((0,\infty),\gamma_\alpha)$ into $\textup{BMO}((0,\infty),\gamma_\alpha)$. 

We can see that $R^\alpha$ is a principal value integral operator. By proceeding as in the proof of \cite[Theorem~1.1]{BFRS} we can see that, for every $f\in L^p((0,\infty),\gamma_\alpha)$, $1\le p<\infty$,
\[
R^\alpha(f)(x)=\lim_{\varepsilon\to 0^+}\int_{|x-y|>\varepsilon,\ y\in (0,\infty)}R^\alpha(x,y)f(y)d\gamma_\alpha(y),\quad \text{a.e. } x\in (0,\infty),
\]
where
\[R^\alpha(x,y)=\frac{1}{\sqrt{\pi}}\int_0^\infty \partial_xW_t^\alpha(x,y)\frac{dt}
{\sqrt{t}},\quad x,y\in (0,\infty),\ x\neq y.
\]
For every $\epsilon>0$, we define the $\epsilon$-truncation of the Riesz transform $R^\alpha$ by
\[
R^\alpha_\epsilon(f)(x)=\int_{|x-y|>\varepsilon,\ y\in (0,\infty)}R^\alpha(x,y)f(y)d\gamma_\alpha(y),\quad x\in (0,\infty).
\]
The maximal Riesz transform $R^\alpha_*$ is  defined by 
\[
R^\alpha_*(f)=\sup_{\epsilon>0}|R^\alpha_\epsilon(f)|.
\]
From the results given by E. Sasso in \cite{SaWeak} we can deduced that the maximal operator~$R^\alpha_*$ is bounded on $L^p((0,\infty),\gamma_\alpha)$, for every $1<p<\infty$, and from $L^1((0,\infty),\gamma_\alpha)$ into $L^{1,\infty}((0,\infty),\gamma_\alpha)$.

We are going to consider the following local maximal Riesz transform operators. For every $a>0$, we define the maximal operator $R^\alpha_{*,a}$ by 
\[R^\alpha_{*,a}(f)(x)= \sup_{0<\epsilon\le am(x)}|R^\alpha_\epsilon(f)(x)|, \quad x\in (0,\infty).\]

Let $\rho>0$. If $\{c_t\}_{t>0}$  is a subset of complex numbers, we define the $\rho$-variation $\mathcal{V}_\rho(\{c_t\}_{t>0})$ of $\{c_t\}_{t>0}$ by
\[
\mathcal{V}_\rho(\{c_t\}_{t>0})=\sup_{0<t_n<t_{n-1}<\dots<t_1,\ n\in\mathbb{N}}\left(\sum_{j=1}^{n-1}|c_{t_j}-c_{t_{j+1}}|^\rho\right)^{1/\rho}.
\]
If $\{T_t\}_{t>0}$ is a family of bounded operators in $L^p((0,\infty),\gamma_\alpha)$, with $1\le p<\infty$, we define the $\rho$-variation operator $\mathcal{V}_\rho(\{T_t\}_{t>0})$ of $\{T_t\}_{t>0}$ by
\[
\mathcal{V}_\rho(\{T_t\}_{t>0})(f)(x)=\mathcal{V}_\rho(\{T_t(f)(x)\}_{t>0}).
\]
Since Bourgain (\cite{Bou}) studied variational inequalities involving martingales (see also \cite{JKRW}), $\rho$-variation operators has been extensively studied in ergodic theory and harmonic analysis. Campbell, Jones, Reinhold and Wierdl (\cite{CJRW1}) proved $L^p$-boundedness properties for $\rho$-variation operators associated to the family of truncations for the Hilbert transform. In \cite{CJRW2} those results were extended by considering Riesz transforms in higher dimensions. In order to obtain $L^p$-boundedness for $\rho$-variation operators it is usual to ask for the condition $\rho>2$ (see \cite{Qi}). For the exponent~$\rho=2$, oscillation operators are commonly considered.

Let $\{t_j\}_{j\in \mathbb{Z}}$ be an increasing sequence of positive real numbers satisfying that $\lim_{j\to -\infty}t_j=0$ and $\lim_{j\to +\infty}t_j=+\infty$. If $\{c_t\}_{t>0}$ is a set of complex numbers, we define the oscillation with respect to $\{t_j\}_{j\in \mathbb{Z}}$ by
\[
\mathcal{O}(\{c_t\}_{t>0},\{t_j\}_{j\in \mathbb{Z}})=\left(\sum_{j=-\infty}^{+\infty}\sup_{t_j\le\epsilon_j<\epsilon_{j+1}<t_{j+1}}|c_{\epsilon_j}-c_{\epsilon_{j+1}}|^2\right)^{1/2}.
\]
If $\{T_t\}_{t>0}$ is a family of bounded operators in $L^p((0,\infty),\gamma_\alpha)$, with $1\le p<\infty$, we define the oscillation operator $\mathcal{O}(\{T_t\}_{t>0},\{t_j\}_{j\in \mathbb{Z}})$ as follows
\[
\mathcal{O}(\{T_t\}_{t>0},\{t_j\}_{j\in \mathbb{Z}})(f)(x)=\mathcal{O}(\{T_t(f)(x)\}_{t>0},\{t_j\}_{j\in \mathbb{Z}}).
\]
$L^p$-boundedness properties of the oscillation operators defined by the family of truncations of Hilbert transform and Euclidean Riesz transforms were established in \cite{CJRW1} and \cite{CJRW2}, respectively.

After \cite{CJRW1} and \cite{CJRW2}, the study of $\rho$-variation and oscillation operators defined by singular integrals has been an active working area (see, for instance, \cite{BCT}, \cite{CDHL}, \cite{CMMTV}, \cite{HMMT}, \cite{MTX1}, \cite{MTX2}, \cite{Mas}, \cite{MasTo} and \cite{MTZ}). Variation and oscillation operators give information about convergence properties for the family $\{T_t\}_{t>0}$.

Being $\{T_t\}_{t>0}$ and $\{t_j\}_{j\in\mathbb{Z}}$ as above, we are going to consider the local $\rho$-variation and oscillation operators defined as follows. Let $a>0$. The $a$-local $\rho$-variation operator $\mathcal{V}_{\rho,a}(\{T_t\}_{t>0})$ is given by
\begin{align*}
    \mathcal{V}_{\rho,a}&(\{T_t\}_{t>0})(f)(x)\\
    &=\sup_{0<t_n<t_{n-1}<\dots<t_1\le am(x),\ n\in\mathbb{N}}\left(\sum_{j=1}^{n-1}|T_{t_j}(f)(x)-T_{t_{j+1}}(f)(x)|^\rho\right)^{1/\rho}.
\end{align*}

The $a$-local oscillation operator $\mathcal{O}_a(\{T_t\}_{t>0},\{t_j\}_{j\in \mathbb{Z}})$ is defined by
\begin{align*}
    \mathcal{O}_a&(\{T_t\}_{t>0},\{t_j\}_{j\in \mathbb{Z}})(f)(x)\\
    &=\left(\sum_{j\in \mathbb{Z},\,\,t_j\le am(x)}\sup_{t_j\le\epsilon_j<\epsilon_{j+1}<t_{j+1}}|T_{\epsilon_j}(f)(x)-T_{\epsilon_{j+1}}(f)(x)|^2\right)^{1/2}.
\end{align*}

Our first result is the following.

\begin{thm}\label{Th1.1}
Let $\alpha>-\frac12$, $a>0$ and $\rho>2$. Suppose that $\{t_j\}_{j\in \mathbb{Z}}$ is an increasing sequence of positive real numbers such that $t_{j+1}\le \theta t_j$, $j\in \mathbb{Z}$, for some $\theta>1$, $\lim_{j\to -\infty}t_j=0$ and $\lim_{j\to +\infty}t_j=+\infty$. The operators $R^\alpha_{*,a}$, $\mathcal{V}_{\rho,a}(\{R^\alpha_\epsilon\}_{\epsilon>0})$, and $\mathcal{O}_a(\{R^\alpha_\epsilon\}_{\epsilon>0},\{t_j\}_{j\in \mathbb{Z}})$ are bounded from $L^\infty((0,\infty),\gamma_\alpha)$ into $\textup{BLO}_a((0,\infty),\gamma_\alpha)$.
\end{thm}

We shall now introduce multiplier operators in the Laguerre setting. A measurable complex function $M$ defined on $[0,\infty)$ is said to be of Laplace transform type when
\[
M(x)=x\int_0^\infty \phi(t)e^{-xt}dt,\quad x>0,\]
where $\phi\in L^\infty(0,\infty)$.

Suppose that $M$ is of Laplace transform type. We denote by $T_M^\alpha$ the spectral multiplier for the Laguerre operator $\Delta_\alpha$ defined by $M-M(0)$. For every ${f\in L^2((0,\infty),\gamma_\alpha)}$, $T_M^\alpha(f)$ is given by
\[
T_M^\alpha(f)=\sum_{k=1}^\infty M(k)c_k^\alpha(f)\mathcal{L}_k^\alpha.
\]

Since $M$ is bounded on $(0,\infty)$, $T_M^\alpha$ is bounded on $L^2((0,\infty),\gamma_\alpha)$. 
Since $\{W_t^\alpha\}_{t>0}$ is a symmetric diffusion semigroup, $T_M^\alpha$ is bounded on $L^p((0,\infty),\gamma_\alpha)$, for every ${1<p<\infty}$ (\cite[Corollary~3,~p.~121]{StLP}). The authors and R. Scotto (\cite[Theorem~1.1~(d)]{BDQS2}) extended the last result establishing variable $L^{p(\cdot)}$-boundedness properties for $T_M^\alpha$. On the other hand, Sasso (\cite{SaSpec}) proved that $T_M^\alpha$ defines a bounded operator from $L^1((0,\infty),\gamma_\alpha)$ into $L^{1,\infty}((0,\infty),\gamma_\alpha)$.  In \cite{BDQS1}, the authors with R. Scotto established the endpoint estimate for $T_M^\alpha$ from $L^\infty((0,\infty),\gamma_\alpha)$ into $\textup{BMO}((0,\infty),\gamma_\alpha)$. 

From \cite[Theorem~1.1]{BFRS} we deduce that there exists a function $\Lambda\in L^\infty(0,\infty)$ such that, for every $f\in L^p((0,\infty),\gamma_\alpha)$, $1\le p<\infty$,
\[
T_M^\alpha(f)(x)=\lim_{\varepsilon\to 0^+}\left(\Lambda(\varepsilon)f(x)+\int_{|x-y|>\varepsilon,\ y\in (0,\infty)}K_\phi^\alpha(x,y)f(y)d\gamma_\alpha(y)\right),\]
for a.e. $x\in (0,\infty)$, where
\[
K_\phi^\alpha(x,y)=-\int_0^\infty\phi(t)\partial_tW_t^\alpha(x,y)dt,\quad x,y\in (0,\infty),\ x\neq y.
\]
A special case of $T_M^\alpha$ is the imaginary power $\Delta_\alpha^{i\alpha}$ that appears when $M_{\eta}(x)=x^{i\eta}$ for $x\in (0,\infty)$ and $\eta\in\mathbb{R}\setminus \{0\}$. For these values of $\eta$,
\[
M_\eta(x)=x\int_0^\infty \phi_\eta(t)e^{-xt}dt,\quad x\in (0,\infty),
\]
where $\phi_\eta(t)=\frac{t^{-i\eta}}{\Gamma(1+i\eta)}$, $t>0$. Note that $|\phi'_\eta(t)|\le C/t$, $t\in (0,\infty)$.

We define, for every $\epsilon>0$, the truncations
\[
Q_{\phi,\epsilon}^\alpha(f)(x)=\int_{|x-y|>\epsilon,\ y\in (0,\infty)}K_\phi^\alpha(x,y)f(y)d\gamma_\alpha(y),\quad x\in (0,\infty),
\]
and consider, for every $a>0$, the $a$-local maximal operator $Q^\alpha_{\phi,*,a}$, which is given by
\[
Q^\alpha_{\phi,*,a}(f)(x)=\sup_{0<\epsilon\le am(x)}|Q_{\phi,\epsilon}^\alpha(f)(x)|.
\]

\begin{thm}\label{Th1.2}
Let $\alpha>-\frac12$ and $a>0$. The maximal operator $Q^\alpha_{\phi,*,a}$ is bounded from $L^\infty((0,\infty),\gamma_\alpha)$ into $\textup{BLO}_a((0,\infty),\gamma_\alpha)$ provided that $|\phi'(t)|\leq C/t$ for some $C>0$, and each $t\in (0,\infty)$.
\end{thm}

The paper is organized as follows. In Section~\ref{sec: BLO} we state the main properties for the spaces $\textup{BLO}_a((0,\infty),\gamma_\alpha)$. In the subsequent sections we prove Theorems \ref{Th1.1} and~\ref{Th1.2}.

Throughout this paper $C$ and $c$ will always denote positive constants than may change in each occurrence.

\section{The spaces \texorpdfstring{$\textup{BLO}_a((0,\infty),\gamma_\alpha)$}{BLOa((0,∞),γα)}}\label{sec: BLO}
In this section we state the main properties of the spaces $\textup{BLO}_a((0,\infty),\gamma_\alpha)$. This properties will be useful in the following sections and they can be proved as the corresponding properties for the Gaussian $\textup{BLO}_a$ space given in \cite[Theorem~3.1, Proposition~3.1 and Theorem~3.2]{LiuY} (see also \cite{Be} for the Euclidean case and \cite{HuYY} for the non-doubling measure case). 

Let $a>0$. The local natural maximal operator $\mathcal{M}^\alpha_a$ associated with the measure~$\gamma_\alpha$ on $(0,\infty)$ is defined by
\[\mathcal{M}^\alpha_a (f)(x) = 
\sup_{I\in \mathcal{B}_a(x)} \frac{1}{\gamma_\alpha(I)} \int_I f(y)d\gamma_\alpha(y),\;x\in(0,\infty),\]
for every measurable function $f$ on $(0,\infty)$ such that $\int_0^\delta |f(y)| d\gamma_\alpha(y)<\infty$, $\delta >0$.

\begin{prop}\label{propo2.1}
    Let $a>0$. There exists $C>0$ such that for every $I\in \mathcal{B}_a$ and every measurable function $f$ on $(0,\infty)$ such that $\|f\|_{*,\alpha}<\infty$,
    \[\frac{1}{\gamma_\alpha(I)}\int_{I} \mathcal{M}^\alpha_a(f)(y) d\gamma_\alpha(y) \leq C  \|f\|_{*,\alpha} + \essinf_{x\in I} \mathcal{M}^\alpha_a(f)(x).\]
    Furthermore, the natural maximal operator $\mathcal{M}^\alpha_a$ defines a bounded operator from $\textup{BMO} ((0,\infty),\gamma_\alpha)$ into $\textup{BLO}_a ((0,\infty),\gamma_\alpha)$.
\end{prop}

The space $\textup{BLO}_a ((0,\infty),\gamma_\alpha)$ can be characterized  by using the local natural maximal operator.

\begin{prop}\label{propo2.2}
    Let $a>0$. A measurable function $f$ belongs to $\textup{BLO}_a ((0,\infty),\gamma_\alpha)$ if and only if $f\in L^1((0,\infty),\gamma_\alpha)$ and $\mathcal{M}^\alpha_a(f)-f\in L^\infty((0,\infty),\gamma_\alpha)$. In addition, we have that 
    \[\|\mathcal{M}^\alpha_a(f)-f\|_{L^\infty((0,\infty),\gamma_\alpha)}
     = \sup_{I\in \mathcal{B}_a}\left(
     \frac{1}{\gamma_\alpha(I)} \int_I f(y) d\gamma_\alpha(y) - \essinf_{x\in I} f(x)\right).\]
\end{prop}

By combining Proposition~\ref{propo2.1} and Proposition~\ref{propo2.2} we can establish the following characterization of $\textup{BLO}_a ((0,\infty),\gamma_\alpha)$ involving the space $\textup{BMO} ((0,\infty),\gamma_\alpha)$ and the local natural maximal operator. 

\begin{prop}
    Let $a>0$. A measurable function $f$ belongs to $\textup{BLO}_a ((0,\infty),\gamma_\alpha)$ if and only if $f=\mathcal{M}^\alpha_a(g) + h$, where $g\in \textup{BMO} ((0,\infty),\gamma_\alpha)$ and $h\in L^{\infty}((0,\infty),\gamma_\alpha)$. Furthermore,
    \[\|f\|_{\textup{BLO}_a ((0,\infty),\gamma_\alpha)} \sim \inf \{ \|g\|_{\textup{BMO} ((0,\infty),\gamma_\alpha)} + \|h\|_{L^\infty((0,\infty),\gamma_\alpha)} \},\]
    where the infimum is taken over all the pairs $(g,h)$ for which $f = \mathcal{M}^\alpha_a(g) + h$ with $(g,h)\in \textup{BMO}((0,\infty),\gamma_\alpha) \times  L^{\infty}((0,\infty),\gamma_\alpha)$.
\end{prop}

\section{Proof of Theorem~\ref{Th1.1}}

\subsection{Local variation operators} \label{subsec: variation}

Let $f\in L^\infty((0,\infty),\gamma_\alpha)$. Since the variation operator $\mathcal{V}_\rho(\{R^\alpha_\epsilon\}_{\epsilon>0})$ is bounded on $L^2((0,\infty),\gamma_\alpha)$ (see~\cite[Theorem~1.3]{BdL}) it follows that
\begin{align*}
        \int_0^\infty  \mathcal{V}_{\rho,a}(\{R^\alpha_\epsilon\}_{\epsilon>0}) (f)(x)d\gamma_\alpha(x) & \leq 
        \left(\int_0^\infty \left(
        \mathcal{V}_\rho( \{R^\alpha_\epsilon\}_{\epsilon>0}) (f)(x)\right)^2d\gamma_\alpha(x)
        \right)^{1/2}
        \\ & 
        \leq C \left(\int_0^\infty |f(x)|^2 d\gamma_\alpha(x) 
        \right)^{1/2}\\
        &\leq C \|f\|_{L^\infty((0,\infty),\gamma_\alpha)}.
    \end{align*}
According to Proposition~\ref{propo2.2} the proof will be finished when we see that 
\begin{equation*}
    \|\mathcal{M}^\alpha_a (\mathcal{V}_{\rho,a}(\{R^\alpha_\epsilon\}_{\epsilon>0}) (f)) - \mathcal{V}_{\rho,a}(\{R^\alpha_\epsilon\}_{\epsilon>0}) (f)\|_{L^{\infty}((0,\infty),\gamma_\alpha)} \leq C \|f\|_{L^{\infty}((0,\infty),\gamma_\alpha)}.
\end{equation*}

Notice that
\begin{align*}
        0 & \leq \mathcal{M}^\alpha_a (\mathcal{V}_{\rho,a}(\{R^\alpha_\epsilon\}_{\epsilon>0}) (f))(x) - \mathcal{V}_{\rho,a}(\{R^\alpha_\epsilon\}_{\epsilon>0}) (f)(x)
        \\ & = \sup_{I \in\mathcal{B}_a(x)} \frac{1}{\gamma_\alpha(I)} \int_{I } \mathcal{V}_{\rho,a}(\{R^\alpha_\epsilon\}_{\epsilon>0}) (f)(z) d\gamma_\alpha(z) - \mathcal{V}_{\rho,a}(\{R^\alpha_\epsilon\}_{\epsilon>0}) (f)(x),
    \end{align*}
for almost every $x\in (0,\infty)$, where $I \in\mathcal{B}_a(x)$ indicates that $I\in \mathcal{B}_a$ and $x\in I$.

Let $x$, $x_0$, $r_0\in (0,\infty)$ such that $I=I(x_0,r_0)\in \mathcal{B}_a(x)$. We decompose $f$ as follows
\begin{equation*}
    f = f \chi_{4I} + f \chi_{(0,\infty)\setminus 4I} = f_1 + f_2.
\end{equation*}

We can write 
\begin{align*}
    \frac{1}{\gamma_\alpha(I)} \int_{I } 
    & 
    \mathcal{V}_{\rho,a}(\{R^\alpha_\epsilon\}_{\epsilon>0}) (f)(z) d\gamma_\alpha(z) - \mathcal{V}_{\rho,a}(\{R^\alpha_\epsilon\}_{\epsilon>0}) (f)(x)
    \\ & 
    \leq \frac{1}{\gamma_\alpha(I)} \int_{I } 
    \mathcal{V}_{\rho,a}(\{R^\alpha_\epsilon\}_{\epsilon>0}) (f_1)(z) d\gamma_\alpha(z)
    \\ & \quad + 
    \frac{1}{\gamma_\alpha(I)} \int_{I } \left(
    \mathcal{V}_{\rho,a}(\{R^\alpha_\epsilon\}_{\epsilon>0}) (f_2)(z)  - \mathcal{V}_{\rho,a}(\{R^\alpha_\epsilon\}_{\epsilon>0}) (f_2)(x) \right) d\gamma_\alpha(z)
    \\ & \quad +
    \mathcal{V}_{\rho,a}(\{R^\alpha_\epsilon\}_{\epsilon>0}) (f_2)(x) - \mathcal{V}_{\rho,a}(\{R^\alpha_\epsilon\}_{\epsilon>0}) (f)(x)
    \\ & := J_1 + J_2 + J_3.
    \end{align*}
By using again that the variation $ \mathcal{V}_\rho(\{R^\alpha_\epsilon\}_{\epsilon>0})$ is bounded on $L^2((0,\infty),\gamma_\alpha)$ we get
\begin{align}
        J_1  & \leq \left(
        \frac{1}{\gamma_\alpha(I)} \int_{I } 
    \left( \mathcal{V}_\rho(\{R^\alpha_\epsilon\}_{\epsilon>0}) (f_1)(z) \right)^2d\gamma_\alpha(z)\right)^{1/2}\nonumber
    \\ & \leq C 
    \left(
        \frac{1}{\gamma_\alpha(I)} \int_{I } 
    \left|  f(z) \right|^2 d\gamma_\alpha(z)\right)^{1/2}
    \leq C \|f\|_{L^\infty((0,\infty),\gamma_\alpha)}.
    \end{align}

Suppose there exists $i_0\in \{1,\dots, n-1\}$ such that $\epsilon_{i_0+1}\leq am(x)<\epsilon_{i_0}$. Thus, for $z\in I$,
\begin{align*}
        &\left(
        \sum_{j=1}^{n-1} |R^\alpha_{\epsilon_{j+1}}(f_2)(z) - R^\alpha_{\epsilon_{j}}(f_2)(z)|^\rho
        \right)^{1/\rho} 
        \\ & \leq 
        \left(
        \sum_{j=1}^{i_0 - 1} |R^\alpha_{\epsilon_{j+1}}(f_2)(z) - R^\alpha_{\epsilon_{j}}(f_2)(z)|^\rho + 
        |R^\alpha_{\epsilon_{i_0}}(f_2)(z) - R^\alpha_{am(x)}(f_2)(z)|^\rho
        \right)^{1/\rho} 
        \\ &\quad + 
        \left(
        \sum_{j=i_0+1}^{n-1} |R^\alpha_{\epsilon_{j+1}}(f_2)(z) - R^\alpha_{\epsilon_{j}}(f_2)(z)|^\rho + 
        |R^\alpha_{\epsilon_{i_0+1}}(f_2)(z) - R^\alpha_{am(x)}(f_2)(z)|^\rho
        \right)^{1/\rho} .
    \end{align*}
Then, recalling that $m(z)\leq Cm(x)$ for every $x,z\in I$, where $C>1$, we obtain
\begin{align*}
         &\mathcal{V}_{\rho,a}(
        \{R^\alpha_\epsilon\}_{\epsilon>0}) (f_2)(z)  - \mathcal{V}_{\rho,a}(\{R^\alpha_\epsilon\}_{\epsilon>0}) (f_2)(x)
        \\ & \leq \sup_{\substack{ 0<\epsilon_n<\dots<\epsilon_1\leq a m(x) \\n\in\mathbb{N}}} \left(
       \sum_{j=1}^{n-1} |R^\alpha_{\epsilon_{j+1}}(f_2)(z) - R^\alpha_{\epsilon_{j}}(f_2)(z)|^\rho   \right)^{1/\rho}
        \\ &\quad + 
        \sup_{\substack{am(x)\leq \epsilon_n<\dots<\epsilon_1<Cam(x)\\n\in\mathbb{N}}}
        \left(
        \sum_{j=1}^{n-1} |R^\alpha_{\epsilon_{j+1}}(f_2)(z) - R^\alpha_{\epsilon_{j}}(f_2)(z)|^\rho
        \right)^{1/\rho}\\
        &\quad - 
      \sup_{\substack{ 0<\epsilon_n<\dots<\epsilon_1\leq a m(x) \\n\in\mathbb{N}}}
        \left(
        \sum_{j=1}^{n-1} |R^\alpha_{\epsilon_{j+1}}(f_2)(x) - R^\alpha_{\epsilon_{j}}(f_2)(x)|^\rho
        \right)^{1/\rho}
        \\ & \leq 
         \sup_{\substack{am(x)\leq \epsilon_n<\dots<\epsilon_1<Cam(x)\\n\in\mathbb{N}}}
        \sum_{j=1}^{n-1} |R^\alpha_{\epsilon_{j+1}}(f_2)(z) - R^\alpha_{\epsilon_{j}}(f_2)(z)|
        \\ & \quad +
\sup_{\substack{0<\epsilon_n<\dots<\epsilon_1\leq am(x)\\n\in\mathbb{N}}} \inf_{\substack{0<\delta_k<\dots<\delta_1\leq am(x)\\k\in\mathbb{N}}}
\left[ \left(\sum_{j=1}^{n-1} |R^\alpha_{\epsilon_{j+1}}(f_2)(z)- R^\alpha_{\epsilon_{j}}(f_2)(z)|^\rho \right)^{1/\rho} \right.
\\ & \qquad - \left.\left(\sum_{j=1}^{k-1} |R^\alpha_{\delta_{j+1}}(f_2)(x) - R^\alpha_{\delta_{j}}(f_2)(x)|^\rho \right)^{1/\rho}\right]
\\ & \leq 
\int_{am(x)<|z-y|<Cam(x)} |R^\alpha(z,y)||f_2(y)| d\gamma_\alpha(y)
\\ & \quad +
\sup_{\substack{0<\epsilon_n<\dots<\epsilon_1\leq am(x)\\n\in\mathbb{N}}} \inf_{\substack{0<\delta_k<\dots<\delta_1\leq am(x)\\k\in\mathbb{N}}}
\left[ \left(\sum_{j=1}^{n-1} |R^\alpha_{\epsilon_{j+1}}(f_2)(z) - R^\alpha_{\epsilon_{j}}(f_2)(z)|^\rho \right)^{1/\rho} \right.
\\ & \qquad - \left.\left(\sum_{j=1}^{k-1} |R^\alpha_{\delta_{j+1}}(f_2)(x) - R^\alpha_{\delta_{j}}(f_2)(x)|^\rho \right)^{1/\rho}\right]\\
& := J_{2,1}(x,z) + J_{2,2}(x,z).
    \end{align*}

If we write 
\begin{equation*}
    R^\alpha(z,y) = e^{\frac{z^2 + y^2}{2}} \mathfrak{R}^\alpha(z,y), \quad z,y\in(0,\infty), \, z\neq y,
\end{equation*}
from \cite[(3.3) and Proposition 3.1]{NS07} we know that
\begin{equation}\label{cota-nucleo-potencia}
    \mathfrak{R}^\alpha(z,y) \leq \frac{C}{\mathfrak{m}_\alpha(I(z,|z-y|))}, \quad z,y\in(0,\infty), \, z\neq y,
\end{equation}
where $d\mathfrak{m}_\alpha(x)=x^{2\alpha+1} dx$. Therefore, for $x,z\in I$, since $m(x)\leq Cm(z)$ and using \cite[(3)]{YY}, we obtain
\begin{align*}
        J_{2,1}(x,z)  & 
        \leq C \int_{am(x)<|z-y|\leq Cam(x)} |f_2(y)| e^{\frac{z^2 - y^2}{2}} \frac{d\mathfrak{m}_\alpha(y)}{\mathfrak{m}_\alpha(I(z,|z-y|))}
        \\ & \leq 
        C \|f\|_{L^\infty((0,\infty),\gamma_\alpha)}\int_{am(x)<|z-y|\leq Cam(x)}  e^{\frac{(z+y)|z-y|}{2}} \frac{d\mathfrak{m}_\alpha(y)}{\mathfrak{m}_\alpha(I(y,|z-y|))}
        \\ & \leq 
        C \|f\|_{L^\infty((0,\infty),\gamma_\alpha)}\int_{am(x)<|z-y|\leq Cam(x)}  e^{Cam(x)(z+Cam(x))} \frac{dy}{|z-y|}
        \\ & \leq 
        C \|f\|_{L^\infty((0,\infty),\gamma_\alpha)}\int_{am(x)<|z-y|\leq Cam(x)}  \frac{dy}{|z-y|}
        \\ & \leq C\|f\|_{L^\infty((0,\infty),\gamma_\alpha)}.
    \end{align*}

On the other hand, for every $x,z\in I$ we have
\begin{align*}
&\left(\sum_{j=1}^{n-1} |R^\alpha_{\epsilon_{j+1}}(f_2)(z) - R^\alpha_{\epsilon_{j}}(f_2)(z)|^\rho \right)^{1/\rho}- \left(\sum_{j=1}^{n-1} |R^\alpha_{\epsilon_{j+1}}(f_2)(x) - R^\alpha_{\epsilon_{j}}(f_2)(x)|^\rho \right)^{1/\rho}
\\ & \leq  
 \left(\sum_{j=1}^{n-1} \left|R^\alpha_{\epsilon_{j+1}}(f_2)(z) - R^\alpha_{\epsilon_{j}}(f_2)(z) -  \left(R^\alpha_{\epsilon_{j+1}}(f_2)(x) - R^\alpha_{\epsilon_{j}}(f_2)(x)\right)\right|^\rho \right)^{1/\rho} 
\\ & \leq 
 \left(\sum_{j=1}^{n-1} \left| \int_{\epsilon_{j+1}<|z-y|<\epsilon_j} (R^\alpha(z,y) - R^\alpha(x,y)) f_2(y) d\gamma_\alpha(y) 
 \right.\right.
\\ & \quad +\left. \left. \left( \int_{\epsilon_{j+1}<|z-y|<\epsilon_j} R^\alpha(x,y) f_2(y) d\gamma_\alpha(y)\right.\right.\right.\\
& \qquad -\left.\left.\left. \int_{\epsilon_{j+1}<|x-y|<\epsilon_j} R^\alpha(x,y) f_2(y) d\gamma_\alpha(y)\right)\right|^\rho 
\right)^{1/\rho} 
\\ & \leq 
\sum_{j=1}^{n-1} \left| \int_{\epsilon_{j+1}<|z-y|<\epsilon_j} (R^\alpha(z,y) - R^\alpha(x,y)) f_2(y) d\gamma_\alpha(y)\right|
\\ & \quad +
\left( \sum_{j=1}^{n-1} \left| \int_0^\infty R^\alpha(x,y) \left(\chi_{\{\epsilon_{j+1}<|x-y|<\epsilon_j\}}(y) \right.\right.\right.
\\ & \qquad -
\left.\left.\left.\chi_{\{\epsilon_{j+1}<|z-y|<\epsilon_j\}}(y)  \right) f_2(y) d\gamma_\alpha(y)\right|^\rho \right)^{1/\rho}.
\end{align*}
Now, by taking supremum, we get
\begin{align*}
J_{2,2}(x,z)&\leq \int_{(0,\infty)\setminus 4I} |R^\alpha(z,y) - R^\alpha(x,y) | |f(y)| d\gamma_\alpha(y) 
\\ & \quad +
\sup_{\substack{0<\epsilon_n<\dots<\epsilon_1\leq am(x)\\n\in\mathbb{N}}} 
\left(\sum_{j=1}^{n-1} \left( \int_0^\infty \left|R^\alpha(x,y)\right| \left|\chi_{\{\epsilon_{j+1}<|x-y|<\epsilon_j\}}(y) \right.\right.\right.
\\ & \qquad -
\left.\left.\left.\chi_{\{\epsilon_{j+1}<|z-y|<\epsilon_j\}}(y)  \right| |f_2(y)| d\gamma_\alpha(y)\right)^\rho \right)^{1/\rho}
\\ & := J_{2,2,1}(x,z) + J_{2,2,2}(x,z).
    \end{align*}

Since (see~\cite[\S~4.3]{BDQS3})
\begin{equation*}
    \sup_{I \in \mathcal{B}_a} \sup_{x,z\in I } \int_{(0,\infty)\setminus 4I} |R^\alpha(z,y) - R^\alpha(x,y) | |f(y)| d\gamma_\alpha(y) < \infty,
\end{equation*}
it follows that
\begin{equation*}
 J_{2,2,1}(x,z) \leq C \|f\|_{L^\infty((0,\infty),\gamma_\alpha)}, \quad x,z\in I .
\end{equation*}

In order to estimate $J_{2,2,2}$ we adapt a procedure developed in~\cite{BFHR}. From \eqref{cota-nucleo-potencia}, for $x,z\in I$,
we obtain
\begin{align*}
    J_{2,2,2}(x,z) & \leq C  \sup_{\substack{0<\epsilon_n<\dots<\epsilon_1\leq am(x)\\n\in\mathbb{N}}} 
\left(\sum_{j=1}^{n-1} \left( \int_0^\infty  \frac{e^{\frac{x^2- y^2}{2}}}{\mathfrak{m}_\alpha(I(y,|x-y|))} \right.\right.
\\ & \qquad \times \left|\chi_{\{\epsilon_{j+1}<|x-y|<\epsilon_j\}}(y)-
\left.\left.\chi_{\{\epsilon_{j+1}<|z-y|<\epsilon_j\}}(y)  \right| |f_2(y)| d\mathfrak{m}_\alpha(y)\right)^\rho \right)^{1/\rho}.
\end{align*}

Let us observe that, if $|x-y|\leq am(x)$, then \[x^2-y^2 \leq |x-y||x+y|\leq am(x)(am(x)+2x)\leq C.\] Also, if $|z-y|\leq am(x)$, then \[|x-y|\leq 2r_0 + am(x) \leq 2am(x_0) + am(x).\] Since $x\in I \in \mathcal{B}_a$, $m(x_0)\leq Cm(x)$ so $|x-y| \leq C am(x)$, and thus $x^2 - y^2 \leq C$, provided that $|z-y| \leq a m(x)$. 

This fact together with~\cite[(3)]{YY} lead to
\begin{align*}
    J_{2,2,2}(x,z) &  \leq \sup_{\substack{0<\epsilon_n<\dots<\epsilon_1\leq am(x)\\n\in\mathbb{N}}} 
\left(\sum_{j=1}^{n-1} \left( \int_0^\infty  \frac{|f_2(y)|}{|x-y|}\left|\chi_{\{\epsilon_{j+1}<|x-y|<\epsilon_j\}}(y) \right.\right.\right.
\\ & \quad -
\left.\left.\left.\chi_{\{\epsilon_{j+1}<|z-y|<\epsilon_j\}}(y)  \right| dy\right)^\rho \right)^{1/\rho}, \quad x,z\in I .
\end{align*}

Let us take $0<\epsilon_n <\dots < \epsilon_1\leq am(x)$ and $j\in\{1,\dots,n-1\}$. Then
\begin{align*}
        \int_0^\infty & \frac{|f_2(y)|}{|x-y|} \left|  \chi_{\{\epsilon_{j+1}<|x-y|<\epsilon_j\}}(y)   - \chi_{\{\epsilon_{j+1}<|z-y|<\epsilon_ j\}}(y) \right|dy 
        \\ & \leq C \left(
        \int_0^\infty \frac{|f_2(y)|}{|x-y|}
        \chi_{\{\epsilon_{j+1}<|x-y|<\epsilon_j\}}(y)\chi_{\{\epsilon_{j+1}<|x-y|<\epsilon_{j+1} + 2r_0\}}(y) dy\right.
        \\ & \quad +
        \int_0^\infty \frac{|f_2(y)|}{|x-y|}
        \chi_{\{\epsilon_{j+1}<|x-y|<\epsilon_j\}}(y)\chi_{\{\epsilon_{j}<|z-y|<\epsilon_{j} + 2r_0\}}(y) dy
        \\ & \quad +
        \int_0^\infty \frac{|f_2(y)|}{|z-y|}
        \chi_{\{\epsilon_{j+1}<|z-y|<\epsilon_j\}}(y)\chi_{\{\epsilon_{j+1}<|z-y|<\epsilon_{j+1} + 2r_0\}}(y) dy
        \\ & \quad + \left.
        \int_0^\infty \frac{|f_2(y)|}{|z-y|}
        \chi_{\{\epsilon_{j+1}<|z-y|<\epsilon_j\}}(y)\chi_{\{\epsilon_{j}<|x-y|<\epsilon_{j} + 2r_0\}}(y) dy \right)
        \\ & = \sum_{l=1}^4 J^{j,l}_{2,2,2} (x,z), \quad x,z\in I .
    \end{align*}

For the above estimate, we have taken into account that, if $\chi_{\{\epsilon_{j+1}<|x-y|\leq \epsilon_j\}} (y) - \chi_{\{\epsilon_{j+1}<|z-y|\leq \epsilon_j\}} (y)\neq 0$, then $\chi_{\{\epsilon_{j+1}<|x-y|\leq \epsilon_j\}} (y)\chi_{\{\epsilon_{j+1}<|z-y|\leq \epsilon_j\}} (y) = 0$, with ${y\in(0,\infty)}$ and $x$, $z\in I $. Since $f_2(y) = 0$ for $y\in 4I$, it follows that $J^{j,l}_{2,2,2} = 0$ when $l=1,3$, $z\in I $ and $r_0\geq \epsilon_{j+1}$. Also, $J^{j,l}_{2,2,2}(x,z) = 0$ when $l=2,4$, $z\in I $ and $r_0\geq \epsilon_j$.

If $z\in I $ and $y\notin 4I$, then $2|x-y|\geq |z-y| \geq \frac12 |x-y|$. H\"older inequality leads to 
\begin{equation*}
        J^{j,l}_{2,2,2} (x,z) \leq C \left(
        \int_0^\infty \chi_{\{\epsilon_{j+1}<|x-y|<\epsilon_j\}}(y) \left(\frac{|f_2(y)|}{|x-y|}\right)^2 dy\right)^{1/2} r_0^{1/2}, \quad z\in I , \, l=1,2;
\end{equation*}
\begin{equation*}
        J^{j,l}_{2,2,2} (x,z) \leq C \left(
        \int_0^\infty \chi_{\{\epsilon_{j+1}<|z-y|<\epsilon_j\}}(y) \left(\frac{|f_2(y)|}{|z-y|}\right)^2 dy\right)^{1/2} r_0^{1/2}, \quad z\in I , \, l=3,4.
\end{equation*}

We obtain
\begin{align*}
        &\left(\sum_{j=1}^{n-1}\left|
        \sum_{l=1}^{4} J^{j,l}_{2,2,2} (x,z)
        \right|^\rho\right)^{1/\rho} \\
        &\quad \leq C  r_0^{1/2} \left[
        \left( \sum_{j=1}^{n-1} 
        \left(
        \int_0^\infty \chi_{\{\epsilon_{j+1}<|x-y|<\epsilon_j\}}(y) \left(\frac{|f_2(y)|}{|x-y|}\right)^2 dy\right)^{\rho/2}\right)^{1/\rho} \right.
        \\ &\qquad + \left. \left( \sum_{j=1}^{n-1}
         \left(
        \int_0^\infty \chi_{\{\epsilon_{j+1}<|z-y|<\epsilon_j\}}(y) \left(\frac{|f_2(y)|}{|z-y|}\right)^2 dy\right)^{\rho/2}
        \right)^{1/\rho} \right]
        \\ & \quad \leq C r_0^{1/2} \|f\|_{L^{\infty}((0,\infty),\gamma_\alpha)}
        \left[
        \left( \sum_{j=1}^{n-1} 
        \int_{(0,\infty)\setminus 4I} \chi_{\{\epsilon_{j+1}<|x-y|<\epsilon_j\}}(y) \frac{dy}{|x-y|^2} \right)^{1/2} \right.
        \\ & \qquad + \left. \left( \sum_{j=1}^{n-1}
        \int_{(0,\infty)\setminus 4I} \chi_{\{\epsilon_{j+1}<|z-y|<\epsilon_j\}}(y) \frac{dy}{|z-y|^2}\right)^{1/2}
         \right]
         \\ & \quad\leq C r_0^{1/2} \|f\|_{L^{\infty}((0,\infty),\gamma_\alpha)}
        \left(  
        \int_{(0,\infty)\setminus 4I}  \frac{dy}{|x-y|^2} +       \int_{(0,\infty)\setminus 4I}  \frac{dy}{|z-y|^2}\right)^{1/2}
        \\ & \quad \leq C  \|f\|_{L^{\infty}((0,\infty),\gamma_\alpha)}.
    \end{align*}
We conclude that $J_{2,2,2}(x,z) \leq C  \|f\|_{L^{\infty}((0,\infty),\gamma_\alpha)}$. By putting together the above estimates we obtain
\begin{equation*}
    \mathcal{V}_{\rho,a} \left(\{R^\alpha_\epsilon\}_{\epsilon>0}\right) (f_2)(z) - 
    \mathcal{V}_{\rho,a} \left(\{R^\alpha_\epsilon\}_{\epsilon>0}\right) (f_2)(x) \leq C  \|f\|_{L^{\infty}((0,\infty),\gamma_\alpha)}, \quad x,z\in I .
\end{equation*}
Here, $C>0$ does not depend on $x$, $z\in I $, so it follows that
\begin{equation*}
    J_2 \leq C \|f\|_{L^{\infty}((0,\infty),\gamma_\alpha)}.
\end{equation*}

We now estimate $J_3$.
Note first that if $|x-y|<r_0$, then $|y-x_0|<2r_0$, so it is clear that 
\begin{equation*}
    \int_{\epsilon_1<|x-y|<\epsilon_2} R^\alpha(x,y)f_2(y) d\gamma_\alpha(y) = 0, \quad 0<\epsilon_1<\epsilon_2\leq r_0.
\end{equation*}

Suppose that $r_0\leq a m(x)$. We have, for any $x\in I$, that
\begin{align*}
        \mathcal{V}_{\rho,a}&\left(\{ R^\alpha_\epsilon\}_{\epsilon>0}\right)(f_2)(x)\\& \leq \sup_{0<\epsilon_n<\dots<\epsilon_1 \leq r_0, \ n\in \mathbb{N}} 
        \left(\sum_{j=1}^{n-1} | R^\alpha_{\epsilon_{j+1}}(f_2)(x) - R^\alpha_{\epsilon_{j}}(f_2)(x)|^\rho\right)^{1/\rho}
        \\ & \quad + \sup_{r_0<\epsilon_n<\dots<\epsilon_1 \leq am(x), \ n\in \mathbb{N}} 
        \left(\sum_{j=1}^{n-1} | R^\alpha_{\epsilon_{j+1}}(f_2)(x) - R^\alpha_{\epsilon_{j}}(f_2)(x)|^\rho\right)^{1/\rho}
        \\ & = \sup_{r_0<\epsilon_n<\dots<\epsilon_1 \leq am(x), \ n\in \mathbb{N}} 
        \left(\sum_{j=1}^{n-1} | R^\alpha_{\epsilon_{j+1}}(f_2)(x) - R^\alpha_{\epsilon_{j}}(f_2)(x)|^\rho\right)^{1/\rho}
        \\ & \leq 
        \mathcal{V}_{\rho,a}\left(\{ R^\alpha_\epsilon\}_{\epsilon>0}\right)(f)(x)
        \\ & \quad + \sup_{r_0<\epsilon_n<\dots<\epsilon_1 \leq am(x), \ n\in \mathbb{N}} 
        \left(\sum_{j=1}^{n-1} | R^\alpha_{\epsilon_{j+1}}(f_1)(x) - R^\alpha_{\epsilon_{j}}(f_1)(x)|^\rho\right)^{1/\rho}.
    \end{align*}

Since, for every $y\in 4I$, $|x-y|\leq 5r_0\leq 5 am(x)$, by using again \eqref{cota-nucleo-potencia} and \cite[(3)]{YY}, we deduce that 
\begin{align*}
\sup_{r_0<\epsilon_n<\dots<\epsilon_1 \leq am(x), \ n\in \mathbb{N}}  &
        \left(\sum_{j=1}^{n-1} | R^\alpha_{\epsilon_{j+1}}(f_1)(x) - R^\alpha_{\epsilon_{j}}(f_1)(x)|^\rho\right)^{1/\rho}
        \\ & \leq C \int_{|x-y|>r_0,\ y\in 4I} \frac{e^{\frac{x^2+y^2}{2}}|f(y)|}{\mathfrak{m}_\alpha(I(y,|x-y|))} d\gamma_\alpha(y)
        \\ & \leq C \int_{|x-y|>r_0,\ y\in 4I}\frac{e^{\frac{x^2-y^2}{2}}|f(y)|}{|x-y|} dy
        \\ & \leq C \|f\|_{L^{\infty}((0,\infty),\gamma_\alpha)}\int_{|x-y|>r_0,\ y\in 4I}\frac{dy}{|x-y|} 
        \\ & \leq C \frac{\|f\|_{L^{\infty}((0,\infty),\gamma_\alpha)}}{r_0}\int_{4I} dy\\
        &  \leq C \|f\|_{L^{\infty}((0,\infty),\gamma_\alpha)},
    \end{align*}
that is, $J_3\leq C \|f\|_{L^{\infty}((0,\infty),\gamma_\alpha)}$ for the case $r_0\leq am(x)$.

When $r_0>am(x)$,
\begin{equation*}
     \mathcal{V}_{\rho,a}\left(\{ R^\alpha_\epsilon\}_{\epsilon>0}\right)(f_2)(x) = 0,
\end{equation*}
so $J_3\leq 0$ in this case.

We conclude that
\begin{equation*}
    J_3 \leq C \|f\|_{L^{\infty}((0,\infty),\gamma_\alpha)}.
\end{equation*}
By combining the above estimates, since the constant $C>0$ does not depend on $x\in (0,\infty)$ or $I\in \mathcal{B}_a(x)$, we get
\begin{equation*}
    \|\mathcal{M}^\alpha_a \left(
    \mathcal{V}_{\rho,a}(\{R^\alpha_\epsilon\}_{\epsilon>0}) (f)\right) - \mathcal{V}_{\rho,a}(\{R^\alpha_\epsilon\}_{\epsilon>0}) (f) \|_{L^{\infty}(0,\infty),\gamma_\alpha} \leq C \|f\|_{L^{\infty}(0,\infty),\gamma_\alpha}. 
\end{equation*}
Thus the proof is finished.

\subsection{Local oscillation operators}

Theorem~\ref{Th1.1} for oscillation operators can be proved by using the procedure developed in the previous section for the variation operator, so we give a sketch of the proof.

According to \cite[Theorem~1.3]{BdL}, the oscillation operator $\mathcal{O}(\{R_\epsilon^\alpha\}_{\epsilon>0}, \{t_j\}_{j\in \mathbb Z})$ is bounded on $L^2((0,\infty),\gamma_\alpha)$. This property implies that $\mathcal{O}_a(\{R_\epsilon^\alpha\}_{\epsilon>0}, \{t_j\}_{j\in \mathbb Z})(f)\in L^1((0,\infty),\gamma_\alpha)$ for every $f\in L^\infty((0,\infty),\gamma_\alpha)$.

In order to prove our result, it is sufficient to find a positive constant $C$ such that, for every $f\in L^\infty((0,\infty),\gamma_\alpha)$,
\begin{align}\label{norm_M-O}
   \nonumber \|\mathcal{M}^\alpha_a (\mathcal{O}_a(\{R^\alpha_\epsilon\}_{\epsilon>0},\{t_j\}_{j\in \mathbb Z}) (f)) - & \mathcal{O}_a(\{R^\alpha_\epsilon\}_{\epsilon>0},\{t_j\}_{j\in \mathbb Z}) (f)\|_{L^{\infty}((0,\infty),\gamma_\alpha)}\\ 
    &\leq C\|f\|_{L^{\infty}((0,\infty),\gamma_\alpha)}.
\end{align}

Fix $f\in L^\infty((0,\infty),\gamma_\alpha)$ and let $x,x_0,r_0\in (0,\infty)$ such that $I=I(x_0,r_0)\in \mathcal{B}_a(x)$. We write $f=f\chi_{4I}+f\chi_{(0,\infty)\setminus 4I}:=f_1+f_2$. 
\begin{align*}
    &\frac{1}{\gamma_\alpha(I)} \int_{I } \mathcal{O}_a(\{R^\alpha_\epsilon\}_{\epsilon>0},\{t_j\}_{j\in \mathbb Z}) (f) (z) d\gamma_\alpha(z) -\mathcal{O}_a(\{R^\alpha_\epsilon\}_{\epsilon>0},\{t_j\}_{j\in \mathbb Z}) (f) (x)\\
    &\leq \frac{1}{\gamma_\alpha(I)} \int_{I } \mathcal{O}_a(\{R^\alpha_\epsilon\}_{\epsilon>0},\{t_j\}_{j\in \mathbb Z}) (f_1) (z) d\gamma_\alpha(z)\\
    &\quad + \frac{1}{\gamma_\alpha(I)} \int_{I } \left[\mathcal{O}_a(\{R^\alpha_\epsilon\}_{\epsilon>0},\{t_j\}_{j\in \mathbb Z}) (f_2) (z)-\mathcal{O}_a(\{R^\alpha_\epsilon\}_{\epsilon>0},\{t_j\}_{j\in \mathbb Z}) (f_2) (x)\right] d\gamma_\alpha(z)\\
    &\quad +\mathcal{O}_a(\{R^\alpha_\epsilon\}_{\epsilon>0},\{t_j\}_{j\in \mathbb Z}) (f_2) (x)-\mathcal{O}_a(\{R^\alpha_\epsilon\}_{\epsilon>0},\{t_j\}_{j\in \mathbb Z}) (f) (x)\\
    &:=J_1+J_2+J_3.
\end{align*}

It is immediate from the $L^2$-boundedness of $\mathcal{O}(\{R_\epsilon^\alpha\}_{\epsilon>0}, \{t_j\}_{j\in \mathbb Z})$ (\cite[Theorem~1.3]{BdL}) that
\[J_1\leq C\|f\|_{L^{\infty}((0,\infty),\gamma_\alpha)}.\]

We now estimate the integrand of $J_2$. For certain $C>1$, we have
\begin{align*}
    \mathcal{O}_a (\{R^\alpha_\epsilon\}_{\epsilon>0}, &\{t_j\}_{j\in \mathbb Z}) (f_2) (z)-\mathcal{O}_a(\{R^\alpha_\epsilon\}_{\epsilon>0},\{t_j\}_{j\in \mathbb Z}) (f_2) (x)\\
    &\leq \left(\sum_{\substack{j\in \mathbb Z\\ t_j\leq C am(x)}} \sup_{t_{j-1}\leq \epsilon_{j-1}<\epsilon_j\leq t_j}\left|R^\alpha_{\epsilon_{j-1}}(f_2)(z)-R^\alpha_{\epsilon_{j}}(f_2)(z)\right|^2 \right)^{1/2}\\
    &\quad -\left(\sum_{\substack{j\in \mathbb Z\\ t_j\leq am(x)}} \sup_{t_{j-1}\leq \epsilon_{j-1}<\epsilon_j\leq t_j}\left|R^\alpha_{\epsilon_{j-1}}(f_2)(x)-R^\alpha_{\epsilon_{j}}(f_2)(x)\right|^2 \right)^{1/2}.
\end{align*}

We define
\begin{equation}\label{j0}
    j_0(x)=\max\{j\in \mathbb Z : t_j\leq am(x)\}
\end{equation}
and also, provided that $t_{j_0(x)+1}\leq Cam(x)$, we consider
\begin{equation}\label{j1}
    j_1(x)=\max\{j\in \mathbb Z : j>j_0(x), t_j\leq Cam(x)\}.
\end{equation}
Thus, when $t_{j_0(x)+1}> Cam(x)$, we can write
\begin{align*}
    \mathcal{O}_a &(\{R^\alpha_\epsilon\}_{\epsilon>0},\{t_j\}_{j\in \mathbb Z}) (f_2) (z)-\mathcal{O}_a(\{R^\alpha_\epsilon\}_{\epsilon>0},\{t_j\}_{j\in \mathbb Z}) (f_2) (x)\\
     &\leq \left(\sum_{\substack{j\in \mathbb Z\\ j\leq j_0(x)}} \sup_{t_{j-1}\leq \epsilon_{j-1}<\epsilon_j\leq t_j}\left|R^\alpha_{\epsilon_{j-1}}(f_2)(z)-R^\alpha_{\epsilon_{j}}(f_2)(z)\right|^2 \right)^{1/2}\\
    &\quad -\left(\sum_{\substack{j\in \mathbb Z\\ j\leq j_0(x)}} \sup_{t_{j-1}\leq \epsilon_{j-1}<\epsilon_j\leq t_j}\left|R^\alpha_{\epsilon_{j-1}}(f_2)(x)-R^\alpha_{\epsilon_{j}}(f_2)(x)\right|^2 \right)^{1/2}\\
    &\leq \left(\sum_{\substack{j\in \mathbb Z\\ j\leq j_0(x)}} \sup_{t_{j-1}\leq \epsilon_{j-1}<\epsilon_j\leq t_j}\left|D(x,z)\right|^2 \right)^{1/2}:=\tilde{J_2}(x,z),
\end{align*}
where \[D(x,z):=R^\alpha_{\epsilon_{j-1}}(f_2)(z)-R^\alpha_{\epsilon_{j}}(f_2)(z)-\left(R^\alpha_{\epsilon_{j-1}}(f_2)(x)-R^\alpha_{\epsilon_{j}}(f_2)(x)\right).\]
On the other hand, if $t_{j_0(x)+1}\leq Cam(x)$, we get
\begin{align*}
     \mathcal{O}_a &(\{R^\alpha_\epsilon\}_{\epsilon>0},\{t_j\}_{j\in \mathbb Z}) (f_2) (z)-\mathcal{O}_a(\{R^\alpha_\epsilon\}_{\epsilon>0},\{t_j\}_{j\in \mathbb Z}) (f_2) (x)\\
     &\leq \tilde{J_2}(x,z)+\left(\sum_{\substack{j\in \mathbb Z\\j_0(x)< j\leq j_1(x)}} \sup_{t_{j-1}\leq \epsilon_{j-1}<\epsilon_j\leq t_j}\left|R^\alpha_{\epsilon_{j-1}}(f_2)(z)-R^\alpha_{\epsilon_{j}}(f_2)(z)\right|^2 \right)^{1/2}\\
     &\leq \tilde{J_2}(x,z)+\int_{\frac a\rho m(x)\leq |z-y|\leq Cam(x)} |R^\alpha(z,y)||f_2(y)| d\gamma_\alpha(y),
\end{align*}
where in the last inequality we have used that $t_{j_0(x)}\leq am(x)\leq t_{j_0(x)+1}\leq \rho t_{j_0(x)}$ with $\rho>1$.

Notice that we can estimate $\tilde{J_2}(x,z)$ in the following form
\begin{align*}
    &\tilde{J_2}(x,z)\\
    &\leq \left(\sum_{\substack{j\in \mathbb Z\\ j\leq j_0(x)}} \sup_{t_{j-1}\leq \epsilon_{j-1}<\epsilon_j\leq t_j}\left|\int_{\epsilon_{j-1}<|z-y|<\epsilon_j} (R^\alpha(z,y)-R^\alpha(x,y))f_2(y) d\gamma_\alpha(y)\right.\right.\\
    &\quad +\left.\left. \int_0^\infty \left(\chi_{\{\epsilon_{j-1}<|z-y|<\epsilon_j\}}(y)-\chi_{\{\epsilon_{j-1}<|x-y|<\epsilon_j\}}(y)\right) R^\alpha(x,y) f_2(y) d\gamma_\alpha(y)\right|^2\right)^{1/2}\\
    &\leq \int_{(0,\infty)\setminus 4I} |R^\alpha(z,y)-R^\alpha(x,y)||f_2(y)| d\gamma_\alpha(y)\\
    &\quad +\left(\sum_{\substack{j\in \mathbb Z\\ j\leq j_0(x)}} \left(\sup_{t_{j-1}\leq \epsilon_{j-1}<\epsilon_j\leq t_j}\int_0^\infty \left|\chi_{\{\epsilon_{j-1}<|z-y|<\epsilon_j\}}(y)-\chi_{\{\epsilon_{j-1}<|x-y|<\epsilon_j\}}(y)\right|\right.\right.\\
    &\qquad \times \left.\left.|R^\alpha(x,y)| |f_2(y)| d\gamma_\alpha(y) \right)^2\right)^{1/2}.
\end{align*}

At this point, we can proceed as in the proof of the corresponding result for variation operators $\mathcal{V}_{\rho,a}$, by using H\"older's inequality with an exponent $s\in (1,2)$ instead of applying it with exponent 2. In this way, we deduce that 
\[J_2\leq C \|f\|_{L^\infty((0,\infty), \gamma_\alpha)}.\]

In order to study $J_3$, we first recall that
\[\int_{\epsilon_1<|x-y|<\epsilon_2} R^\alpha(x,y) f_2(y)d\gamma_\alpha(y), \quad 0<\epsilon_1<\epsilon_2\leq r_0.\]
Then, if $r_0\geq am(x)$, we obtain
\begin{equation}\label{osc-zero}
    \mathcal{O}_a(\{R^\alpha_\epsilon\}_{\epsilon>0},\{t_j\}_{j\in \mathbb Z}) (f_2) (x)=0.
\end{equation}

Suppose now that $r_0<am(x)$ and define $j_0(x)$ as in \eqref{j0}. If $t_{j_0(x)}\leq r_0$, we again have \eqref{osc-zero}. If not, we define $j_1=\max\{j\in \mathbb Z: t_{j_1}\leq r_0\}$. Then
\begin{align*}
    \mathcal{O}_a (\{R^\alpha_\epsilon\}_{\epsilon>0},& \{t_j\}_{j\in \mathbb Z}) (f_2) (x)\\
    &=\left(\sum_{j=j_1+1}^{j_0(x)} \sup_{t_{j-1}\leq \epsilon_{j-1}<\epsilon_j\leq t_j} |R^\alpha_{\epsilon_{j-1}}(f_2)(x)-R^\alpha_{\epsilon_{j}}(f_2)(x)|^2\right)^{1/2}\\
    &\leq \mathcal{O}_a(\{R^\alpha_\epsilon\}_{\epsilon>0},\{t_j\}_{j\in \mathbb Z}) (f) (x)\\
    &\quad +\left(\sum_{j=j_1+1}^{j_0(x)} \sup_{t_{j-1}\leq \epsilon_{j-1}<\epsilon_j\leq t_j} |R^\alpha_{\epsilon_{j-1}}(f_1)(x)-R^\alpha_{\epsilon_{j}}(f_1)(x)|^2\right)^{1/2}.
\end{align*}
Since $t_{j_1}\leq r_0\leq t_{j+1}\leq \rho t_{j_1}$, it follows that
\begin{align*}
    &\left(\sum_{j=j_1+1}^{j_0(x)} \sup_{t_{j-1}\leq \epsilon_{j-1}<\epsilon_j\leq t_j} |R^\alpha_{\epsilon_{j-1}}(f_1)(x)-R^\alpha_{\epsilon_{j}}(f_1)(x)|^2\right)^{1/2}\\
    &\quad \leq \sum_{j=j_1+1}^{j_0(x)} \sup_{t_{j-1}\leq \epsilon_{j-1}<\epsilon_j\leq t_j} |R^\alpha_{\epsilon_{j-1}}(f_1)(x)-R^\alpha_{\epsilon_{j}}(f_1)(x)|\\
    &\quad \leq C \int_{|x-y|>r_0/\rho} \frac{e^{\frac{x^2+y^2}{2}}|f(y)|}{\mathfrak{m}_\alpha (I(y,|x-y|))} d\gamma_\alpha(y)\\
    &\quad\leq C\|f\|_{L^{\infty}((0,\infty),\gamma_\alpha)},
\end{align*}
where we have used again the bound given in \eqref{cota-nucleo-potencia} and \cite[(3)]{YY}.

We conclude that 
\[J_3\leq C\|f\|_{L^{\infty}((0,\infty),\gamma_\alpha)}.\]

By putting together all of the above estimates, we get \eqref{norm_M-O} and the proof of Theorem~\ref{Th1.1} for local oscillation operators is completed.

\subsection{Local maximal Riesz transform}

We firstly prove that $R^\alpha_{*,a}$ is bounded on $L^p((0,\infty),\gamma_\alpha)$ for every $1<p<\infty$. In order to do so, we need to decompose, for every $\epsilon>0$, the truncated integral $R^\alpha_\epsilon$ into two parts, called local and global parts (see \cite{SaWeak}).

For every $\tau>0$, we consider the sets
\[L_\tau=\left\{(x,y,s)\in \zinfa : \sqrt{q_-(x,y,s)}\leq \frac{a(1+\alpha)\tau}{1+x+y}\right\}\]
and
\[G_\tau=(\zinfa )\setminus L_\tau.\]
Here, and in the sequel, we denote $q_{\pm}(x,y,s)=x^2+y^2\pm 2xys$, for $x,y\in (0,\infty)$ and $s\in (-1,1)$. 

We choose a function $\varphi\in C^\infty (\zinfa)$ such that $0\leq \varphi\leq 1$,
\[\varphi(x,y,s)=\begin{cases}1, & (x,y,s)\in L_1,\\ 0, & (x,y,s)\in G_2, \end{cases}\]
and
\[|\partial_x \varphi(x,y,s)|+|\partial_y \varphi(x,y,s)|\leq \frac{C}{\sqrt{q_-(x,y,s)}}, \quad  x,y\in (0,\infty), s\in (-1,1).\]

We define, for each $\epsilon>0$ and $x\in (0,\infty)$
\begin{align*}
    R^{\alpha,\loc}_\epsilon (f)(x)&=\int_{|x-y|>\varepsilon,\ y\in (0,\infty)}R^{\alpha,\loc}(x,y)f(y)d\gamma_\alpha(y),\\ R^{\alpha,\glob}_\epsilon (f)(x)&=R^\alpha_\epsilon (f)(x)-R^{\alpha,\loc}_\epsilon (f)(x),\
\end{align*}
where
\[R^{\alpha,\loc} (x,y)=\int_{-1}^1 R^{\alpha,\loc} (x,y,s) \Pi_\alpha(s)ds,\quad x,y\in (0,\infty)\]
and, for $x,y\in (0,\infty)$, $s\in (-1,1)$,
\[R^{\alpha,\loc} (x,y,s)=-\frac{2}{\sqrt{\pi}}\int_0^\infty \frac{e^{-t(\alpha+2)}\left(e^{-t}x-ys\right)}{(1-e^{-2t})^{\alpha+2}} e^{-\frac{q_-(e^{-t}x,y,s)}{1-e^{-2t}}+y^2} \varphi(x,y,s) \frac{dt}{\sqrt{t}}.\]

We also consider the maximal operators associated with the above,
\[R^{\alpha,\loc}_* (f)=\sup_{\epsilon>0} \left|R^{\alpha,\loc}_\epsilon (f)\right|, \quad R^{\alpha,\glob}_* (f)=\sup_{\epsilon>0} \left|R^{\alpha,\glob}_\epsilon (f)\right|,\]
which clearly verify
\[R^\alpha_* (f)\leq R^{\alpha,\loc}_* (f)+R^{\alpha,\glob}_* (f).\]

According to \cite[\S~3.1]{BdL} (see also \cite[Proposition~3.1]{SaWeak}), we have that
\[R^{\alpha,\glob}_* (f)(x)\leq C\int_0^\infty  K^\alpha(x,y)f(y) d\gamma_\alpha(y), \quad x\in (0,\infty),\]
where
\[K^\alpha(x,y)=\int_{-1}^1 K^\alpha(x,y,s) \chi_{G_1}(x,y,s) \Pi_\alpha(s) ds, \quad x,y\in (0,\infty),\]
and, for $x,y\in (0,\infty)$ and $s\in (-1,1)$,
\begin{equation}\label{Kalfa}
    K^\alpha(x,y,s)=\begin{cases}
1, & s<0, \\
\left(\frac{q_+(x,y,s)}{q_-(x,y,s)}\right)^{\frac{\alpha+1}{2}}\exp\left(\frac{x^2+y^2-\sqrt{q_-(x,y,s)q_+(x,y,s)}}{2}\right), & s\geq 0.
\end{cases}
\end{equation}
It follows that $R^{\alpha,\glob}_*$ is bounded on $L^p((0,\infty),\gamma_\alpha)$ for every $1<p<\infty$ (see \cite[\S~3.1]{BdL}).

We recall that the measure $\mathfrak{m}_\alpha$ defined in Section~\ref{subsec: variation} has the doubling property on $(0,\infty)$. Therefore, by \cite[(18)~and~(19)]{BdL}, $e^{-y^2}R_{\alpha}^\loc(x,y)$, for $x,y\in (0,\infty)$, is an $\mathfrak{m}_\alpha$-standard Calderón-Zygmund kernel, that is, for every $x,y\in (0,\infty)$, $x\neq y$,
\[\left|e^{-y^2}R_{\alpha}^\loc(x,y)\right|\leq \frac{C}{\mathfrak{m}_\alpha(I(x,|x-y|))},\]
and
\[\left|\partial_x\left[e^{-y^2}R_{\alpha}^\loc(x,y)\right]\right|+\left|\partial_y\left[e^{-y^2}R_{\alpha}^\loc(x,y)\right]\right|\leq \frac{C}{|x-y|\mathfrak{m}_\alpha(I(x,|x-y|))}.\]

If we define the operators $R^{\alpha,\loc}$ and $R^{\alpha,\glob}$ in the obvious way, we can see as above that the later is bounded on $L^p((0,\infty),\gamma_\alpha)$ for every $1<p<\infty$. Since  $R^\alpha$ is also bounded on $L^p((0,\infty),\gamma_\alpha)$ for every $1<p<\infty$ (see \cite[Theorem~13]{No}), we conclude that $R^{\alpha,\loc}$ is bounded on $L^p((0,\infty),\gamma_\alpha)$ for every $1<p<\infty$. By proceeding as in \cite[\S~2]{BDQS1}, we deduce that $R^{\alpha,\loc}$ is bounded on $L^p((0,\infty),\mathfrak{m}_\alpha)$ for every $1<p<\infty$. Moreover, since $R^{\alpha,\loc}$ is an $\mathfrak{m}_\alpha$-Calderón-Zygmund operator, $R^{\alpha,\loc}_*$ is bounded on $L^p((0,\infty),\mathfrak{m}_\alpha)$ for every $1<p<\infty$. By using again the arguments given in \cite[\S~2]{BDQS1}, we get that $R^{\alpha,\loc}_*$ is bounded on $L^p((0,\infty),\gamma_\alpha)$ for every $1<p<\infty$.

It follows now that $R^\alpha_*$ is  bounded on $L^p((0,\infty),\gamma_\alpha)$ for every $1<p<\infty$. Particularly, using this property for $p=2$, for any $f\in L^\infty((0,\infty),\gamma_\alpha)$, 
\[\|R^\alpha_{*,a}\|_{L^1((0,\infty),\gamma_\alpha)}\leq C\|f\|_{L^\infty((0,\infty),\gamma_\alpha)}.\]

We recall that, from \eqref{cota-nucleo-potencia},
\[|R^\alpha(x,y)|\leq C\frac{e^{\frac{x^2+y^2}{2}}}{\mathfrak{m}_\alpha(I(x,|x-y|))}, \quad x,y\in (0,\infty),\ x\neq y,\]
and also we can see (as in \cite[\S~4.3]{BDQS3}) that
\[\sup_{I \in \mathcal{B}_a} \sup_{x\in I } r_0 \int_{(0,\infty)\setminus 2I} |\partial_x R^\alpha(x,y)| d\gamma_\alpha(y)<\infty.\]
By proceeding as in the proof of \cite[Theorem~4.1]{LiuY}, it yields
\[\sup_{I \in \mathcal{B}_a}\left\|\mathcal{M}^\alpha_a\left(R^\alpha_{*,a}(f)\right)-R^\alpha_{*,a}(f)\right\|_{L^\infty((0,\infty),\gamma_\alpha)}\leq C\|f\|_{L^\infty((0,\infty),\gamma_\alpha)},\]
meaning that $R^\alpha_{*,a}$ is bounded from $L^\infty((0,\infty),\gamma_\alpha)$ into $\textup{BLO}_a((0,\infty),\gamma_\alpha)$.

\section{Proof of Theorem~\ref{Th1.2}}

In this section, we will study $L^\infty((0,\infty),\gamma_\alpha)$-$\textup{BLO}_a((0,\infty),\gamma_\alpha)$ estimates for the $a$-local maximal operator
\begin{align*}
    Q^\alpha_{\phi,*,a}(f)(x)&=\sup_{0<\epsilon\le am(x)}|Q_{\phi,\epsilon}^\alpha(f)(x)|\\
    &=\sup_{0<\epsilon\le am(x)}\left|\int_{|x-y|>\epsilon, \ y\in (0,\infty)}K_\phi^\alpha(x,y)f(y)d\gamma_\alpha(y)\right|,
\end{align*}
for $x\in (0,\infty)$ and $a>0$.

We recall that
\[K_\phi^\alpha(x,y)=-\int_0^\infty \phi(t) \partial_t W_t^\alpha(x,y) dy, \quad x,t,\in (0,\infty),\ x\neq y,\]
being
\[W_t^\alpha(x,y)=\left(\frac{e^{-t}}{1-e^{-2t}}\right)^{\alpha+1} \int_{-1}^1 e^{-\frac{q_-(e^{-t}x,y,s)}{1-e^{-2t}}+y^2}\Pi_\alpha(s) ds, \quad x,y,t\in (0,\infty).\]

Firstly, we shall see that $Q^\alpha_{\phi,*,a}$ is bounded on $L^p((0,\infty),\gamma_\alpha)$ for every~${1<p<\infty}$. We define, for $x,y,t\in (0,\infty)$,
\[W_t^{\alpha,\loc}(x,y)=\left(\frac{e^{-t}}{1-e^{-2t}}\right)^{\alpha+1} \int_{-1}^1 e^{-\frac{q_-(e^{-t}x,y,s)}{1-e^{-2t}}+y^2}\varphi(x,y,s)\Pi_\alpha(s) ds,\]
and
\[W_t^{\alpha,\glob}(x,y)=W_t^\alpha(x,y)-W_t^{\alpha,\loc}(x,y).\]
In terms of these, we consider $K_\phi^{\alpha,\loc}$ and $K_\phi^{\alpha,\glob}$ given as $K_\phi^\alpha$ but with $W_t^\alpha$ replaced by $W_t^{\alpha,\loc}$ and $W_t^{\alpha,\glob}$, respectively. Similarly, we define $Q_{\phi,*,a}^{\alpha,\loc}$ and $Q_{\phi,*,a}^{\alpha,\glob}$ by putting $K_\phi^{\alpha,\loc}$ and $K_\phi^{\alpha,\glob}$ instead of $K_\phi^\alpha$, respectively.

We will first deal with $Q_{\phi,*,a}^{\alpha,\glob}$. Notice that, for $x,y,t\in (0,\infty)$ and $s\in(-1,1)$
\begin{align*}
    \partial_t&\left[\left(\frac{e^{-t}}{1-e^{-2t}}\right)^{\alpha+1} \exp\left(-\frac{q_-(e^{-t}x,y,s)}{1-e^{-2t}}\right)\right]\\
    &=P_{x,y,s}\left(e^{-t}\right) \left(\frac{e^{-t}}{1-e^{-2t}}\right)^{\alpha+1}\exp\left(-\frac{q_-(e^{-t}x,y,s)}{1-e^{-2t}}\right), 
    \end{align*}
where, for every $x,y\in (0,\infty)$ and $s\in (-1,1)$, $P_{x,y,s}$ is a polynomial whose degree is at most four. Hence,
\begin{align*}
    |K_\phi^{\alpha,\glob}(x,y)|&\leq C \int_{-1}^1 \sup_{t>0} \left(\frac{e^{-t}}{1-e^{-2t}}\right)^{\alpha+1} e^{-\frac{q_-(e^{-t}x,y,s)}{1-e^{-2t}}+y^2}\chi_{L_1^c}(x,y,s)\Pi_\alpha(s) ds\\
    &\leq C \int_{-1}^1 K^\alpha(x,y,s) \chi_{L_1^c}(x,y,s)\Pi_\alpha(s) ds, \quad x,y\in (0,\infty),
\end{align*}
being $K^\alpha(x,y,s)$ as in \eqref{Kalfa}, for $(x,y,s)\in L_1^c$. 

From \cite[\S~3.1]{BdL}, it follows that the operator whose kernel is the one on the right-hand side is bounded on $L^p((0,\infty),\gamma_\alpha)$ for every~${1<p<\infty}$, and so will be~$Q_{\phi,*,a}^{\alpha,\glob}$.

Furthermore, for every $f\in L^p((0,\infty),\gamma_\alpha)$, $1<p<\infty$,
\[\lim_{\varepsilon\to 0^+} \int_{|x-y|>\varepsilon,\ y\in (0,\infty)} K_\phi^{\alpha,\glob}(x,y) f(y) d\gamma_\alpha(y)=\int_0^\infty K_\phi^{\alpha,\glob}(x,y) f(y) d\gamma_\alpha(y),\]
for a.e. $x\in (0,\infty)$.

We now consider the operators 
\[
T_M^{\alpha,\loc}(f)(x)=\lim_{\varepsilon\to 0^+}\left(\Lambda(\varepsilon)f(x)+\int_{|x-y|>\varepsilon,\ y\in (0,\infty)}K_\phi^{\alpha,\loc}(x,y)f(y)d\gamma_\alpha(y)\right),\]
and
\[
T_M^{\alpha,\glob}(f)(x)=\int_0^\infty K_\phi^{\alpha,\glob}(x,y) f(y) d\gamma_\alpha(y),\]
for a.e. $x\in (0,\infty)$.
Since $T_M^\alpha$ and $T_M^{\alpha,\glob}$ are both bounded on $L^2((0,\infty),\gamma_\alpha)$ (\cite[Proposition~3]{SaSpec}), also $T_M^{\alpha,\loc}$ is bounded on $L^2((0,\infty),\gamma_\alpha)$. Moreover, for every $f\in L^\infty((0,\infty),\gamma_\alpha)$
\[
T_M^{\alpha,\loc}(f)(x)=\int_0^\infty K_\phi^{\alpha,\loc}(x,y)f(y)d\gamma_\alpha(y),\quad x\notin \supp(f).\]

Let us now consider $\mathbb{K}_\phi^\alpha(x,y):=e^{-y^2} K_\phi^{\alpha,\loc}(x,y)$, for $x,y\in (0,\infty)$. We have
\begin{align*}
    \mathbb{K}_\phi^\alpha(x,y)&=(\alpha+1)\int_0^\infty \varphi(t) \left(\frac{e^{-t}}{1-e^{-2t}}\right)^{\alpha+1} \int_{-1}^1 e^{-\frac{q_-(e^{-t}x,y,s)}{1-e^{-2t}}}\varphi(x,y,s) \Pi_\alpha(s) ds dt\\
    &\quad -\int_0^\infty \varphi(t) e^{-t(\alpha+1)} \int_{-1}^1 \partial_t \left[\frac{e^{-\frac{q_-(e^{-t}x,y,s)}{1-e^{-2t}}}}{(1-e^{-2t})^{\alpha+1}}\right]\varphi(x,y,s)\Pi_\alpha(s) ds dt\\
    &:=\mathbb{K}_{\phi,1}^\alpha(x,y)+\mathbb{K}_{\phi,2}^\alpha(x,y), \quad x,y\in (0,\infty).
\end{align*}
As in \cite[\S~7]{BDQS1}, we can prove that
\[|\mathbb{K}_{\phi,2}^\alpha(x,y)|\leq \frac{C}{\mathfrak{m}_\alpha(I(x,|x-y|))}, \quad x,y\in (0,\infty),\ x\neq y,\]
and
\[|\partial_x\mathbb{K}_{\phi,2}^\alpha(x,y)|+|\partial_y\mathbb{K}_{\phi,2}^\alpha(x,y)|\leq \frac{C}{|x-y|\mathfrak{m}_\alpha(I(x,|x-y|))}, \quad x,y\in (0,\infty),\ x\neq y.\]

On the other hand, using \cite[(2.6)]{SaSpec}, i.e., $q_-(e^{-t}x,y,s)\geq q_-(x,y,s)-2(1-e^{-2t})$, for every $(x,y,s)\in N_1$, and the estimates obtained in \cite[p.~12~and~Lemma~3.1]{BCN}, we get
\begin{enumerate}
    \item 
    \begin{align*}
        |\mathbb{K}_{\phi,1}^\alpha(x,y)|&\leq C \int_0^\infty |\varphi(t)| \left(\frac{e^{-t}}{1-e^{-2t}}\right)^{\alpha+1} \int_{-1}^1 e^{-\frac{q_-(x,y,s)}{1-e^{-2t}}} \Pi_\alpha(s) dsdt\\
        & \leq C \int_0^\infty |\varphi(t)| e^{-t(\alpha+1)} dt\int_{-1}^1 \frac{\Pi_\alpha(s)}{q_-(x,y,s)^{\alpha+1}} ds\\
        &\leq \frac{C}{\mathfrak{m}_\alpha(I(x,|x-y|))}, \quad x,y\in (0,\infty),\ x\neq y;
    \end{align*}
    \item by \cite[Lemma~3.4~(E8)]{BDQS3},
    \begin{align*}
        |\partial_x \mathbb{K}_{\phi,1}^\alpha(x,y)|&\leq C \int_0^\infty |\varphi(t)| \left(\frac{e^{-t}}{1-e^{-2t}}\right)^{\alpha+1} \int_{-1}^1 e^{-\frac{q_-(e^{-t}x,y,s)}{1-e^{-2t}}}\\
        &\qquad \times \left[\frac{e^{-t}|e^{-t}x-ys|}{1-e^{-2t}}\varphi(x,y,s)+|\partial_x \varphi(x,y,s)|\right] \Pi_\alpha(s) dsdt\\
        &\leq C \int_0^\infty |\varphi(t)| \left(\frac{e^{-t}}{1-e^{-2t}}\right)^{\alpha+1} e^{-t}\int_{-1}^1 e^{-\frac{q_-(e^{-t}x,y,s)}{1-e^{-2t}}}\\
        &\qquad \times \left[\frac{\sqrt{q_-(e^{-t}x,y,s)}}{1-e^{-2t}}+\frac{1}{\sqrt{q_-(x,y,s)}}\right] \Pi_\alpha(s) dsdt\\
        &\leq C \int_0^\infty |\varphi(t)| e^{-t(\alpha+2)} \int_{-1}^1 \left[\frac{e^{-\frac{q_-(x,y,s)}{2(1-e^{-2t})}}}{(1-e^{-2t})^{\alpha+3/2}}\right.\\
        &\qquad +\left.\frac{e^{-\frac{q_-(x,y,s)}{1-e^{-2t}}}}{(1-e^{-2t})^{\alpha+1} \sqrt{q_-(x,y,s)}}\right]\Pi_\alpha(s) dsdt\\
        &\leq C \int_0^\infty |\varphi(t)| e^{-t(\alpha+2)} dt\int_{-1}^1 \frac{\Pi_\alpha(s)}{q_-(x,y,s)^{\alpha+3/2}} ds\\
        &\leq \frac{C}{|x-y|\mathfrak{m}_\alpha(I(x,|x-y|))}, \quad x,y\in (0,\infty),\ x\neq y;
    \end{align*}
    \item by \cite[Lemma~3.4~(E7)]{BDQS3} (with $x$ and $y$ interchanged), and proceeding like before,
    \begin{align*}
        |\partial_y \mathbb{K}_{\phi,1}^\alpha(x,y)|&\leq C \int_0^\infty |\varphi(t)| \left(\frac{e^{-t}}{1-e^{-2t}}\right)^{\alpha+1} \int_{-1}^1 e^{-\frac{q_-(e^{-t}x,y,s)}{1-e^{-2t}}}\\
        &\qquad \times \left[\frac{|y-e^{-t}xs|}{1-e^{-2t}}\varphi(x,y,s)+|\partial_y \varphi(x,y,s)|\right] \Pi_\alpha(s) ds dt\\
        &\leq C \int_0^\infty |\varphi(t)| \left(\frac{e^{-t}}{1-e^{-2t}}\right)^{\alpha+1} \int_{-1}^1 e^{-\frac{q_-(e^{-t}x,y,s)}{1-e^{-2t}}}\\
        &\qquad \times \left[\frac{\sqrt{q_-(e^{-t}x,y,s)}}{1-e^{-2t}}+\frac{1}{\sqrt{q_-(x,y,s)}}\right] \Pi_\alpha(s) ds dt\\
        &\leq C \int_0^\infty |\varphi(t)| e^{-t(\alpha+1)} dt\int_{-1}^1 \frac{\Pi_\alpha(s)}{q_-(x,y,s)^{\alpha+3/2}} ds \\
        &\leq \frac{C}{|x-y|\mathfrak{m}_\alpha(I(x,|x-y|))}, \quad x,y\in (0,\infty),\ x\neq y.
    \end{align*}
\end{enumerate}

All of the above proves that $T_M^{\alpha,\loc}$ is an $\mathfrak{m}_\alpha$-Calderón-Zygmund operator. Therefore, $T_M^{\alpha,\loc}$ is bounded on $L^p((0,\infty),\mathfrak{m}_\alpha)$ for every $1<p<\infty$, which yields $Q_{\phi,*,a}^{\alpha,\loc}$ is also bounded on $L^p((0,\infty),\mathfrak{m}_\alpha)$ for every~${1<p<\infty}$. The arguments in \cite[\S~2]{BDQS1} allow us to deduce that $Q_{\phi,*,a}^{\alpha,\loc}$ is also bounded on $L^p((0,\infty),\gamma_\alpha)$ for every~${1<p<\infty}$.

Finally, we conclude that $Q^\alpha_{\phi,*,a}$ is bounded on $L^p((0,\infty),\gamma_\alpha)$ for any $1<p<\infty$.

\begin{rem} We can also prove that $Q^\alpha_{\phi,*,a}$ is bounded from $L^1((0,\infty),\gamma_\alpha)$ to $L^{1,\infty}((0,\infty),\gamma_\alpha)$. Actually, it is sufficient at this moment to know that $Q^\alpha_{\phi,*,a}$ is bounded on $L^{p_0}((0,\infty),\gamma_\alpha)$ for some $1<p_0<\infty$.
\end{rem}

We have proved above that
\[|K_\phi^{\alpha,\loc}(x,y)|\leq C \frac{e^{y^2}}{\mathfrak{m}_\alpha(I(x,|x-y|))}, \quad x,y\in (0,\infty),\ x\neq y.\]
We also saw that
\[|K_\phi^{\alpha,\glob}(x,y)|\leq C \int_{-1}^1 K^\alpha(x,y,s) \chi_{L_1^c}(x,y,s)\Pi_\alpha(s) ds, \quad x,y\in (0,\infty),\]
where $K^\alpha(x,y,s)$ was defined in \eqref{Kalfa}. 
It is easy to see that, for any $(x,y,s)\in L_1^c$,
\[|K^\alpha(x,y,s)|\leq C\begin{cases}
1, & s\in (-1,0), \\
\frac{\exp\left(\frac{x^2+y^2}{2}\right)}{q_-(x,y,s)^{\alpha+1}}, & s\in [0,1).
\end{cases}\]

Moreover, for any fixed constant $c>0$, if $x,y\in (0,\infty)$ with $|x-y|\leq cam(x)$ and $s\in (-1,1)$,
\begin{align*}
    q_-(x,y,s)&=(x-y)^2+2xy(1-s)\leq (x-y)^2+4y(|x-y|+y)\\
    &\leq 5(x-y)^2+4y^2+4|x-y|x\leq C(1+y^2),
\end{align*}
which yields
\[|K^\alpha(x,y,s)|\leq \frac{C}{q_-(x,y,s)^{\alpha+1}}\begin{cases}
(1+y^2)^{\alpha+1}, & s\in (-1,0), \\
\exp\left(\frac{x^2+y^2}{2}\right), & s\in [0,1),
\end{cases}\]
for any $(x,y,s)\in L_1^c$ with $|x-y|\leq cam(x)$.

According to \cite[Lemma~3.1]{BCN}, we obtain
\[|K_\phi^{\alpha,\glob}(x,y)|\leq C \frac{e^{\frac{x^2+y^2}{2}}}{\mathfrak{m}_\alpha(I(x,|x-y|))}, \quad x,y\in (0,\infty), 0<|x-y|\leq cam(x).\]
Hence, we conclude that
\[|K_\phi^\alpha(x,y)|\leq C \frac{e^{\frac{x^2+y^2}{2}}}{\mathfrak{m}_\alpha(I(x,|x-y|))}, \quad x,y\in (0,\infty), 0<|x-y|\leq cam(x).\]

We are going to prove now that
\[\sup_{I \in \mathcal{B}_1} \sup_{x\in I } \int_{(0,\infty)\setminus 2I} |\partial_x K_\phi^\alpha(x,y)|d\gamma_\alpha(y) <\infty.\]

By partial integration, we have that
\[K_\phi^\alpha(x,y)=\int_0^\infty \phi'(t) W_t^\alpha(x,y) dt, \quad x,y,\in (0,\infty),\]
and thus,
\begin{align*}
    \partial_x K_\phi^\alpha(x,y) & =\int_0^\infty \phi'(t) \partial_x W_t^\alpha(x,y) dt\\
    &=-2\int_0^\infty \phi'(t) \left(\frac{e^{-t}}{1-e^{-2t}}\right)^{\alpha+1}\\
    &\quad \times \int_{-1}^1 e^{-\frac{q_-(e^{-t}x,y,s)}{1-e^{-2t}}+y^2} \frac{e^{-t}(e^{-t}x-ys)}{1-e^{-2t}} \Pi_\alpha(s) ds dt, \quad x,y\in (0,\infty).
\end{align*}
By \cite[Lemma~3.4~(E8)]{BDQS3}, $|e^{-t}x-ys|\leq \sqrt{q_-(e^{-t}x,y,s)}$ for every $x,y\in (0,\infty)$ and $s\in (-1,1)$. Then, using the hypothesis on $\phi'$,
\begin{align*}
    |\partial_x K_\phi^\alpha(x,y)|&\leq C \int_0^\infty |\phi'(t)| \frac{e^{-t(\alpha+2)}}{(1-e^{-2t})^{\alpha+1}} \\
    &\quad \times\int_{-1}^1  e^{-\frac{q_-(e^{-t}x,y,s)}{1-e^{-2t}}+y^2} \sqrt{q_-(e^{-t}x,y,s)}\Pi_\alpha(s) ds dt\\
    &\leq C\int_0^\infty \frac1t \frac{e^{-t(\alpha+2)}}{(1-e^{-2t})^{\alpha+3/2}}\int_{-1}^1  e^{-c\frac{q_-(e^{-t}x,y,s)}{1-e^{-2t}}+y^2}\Pi_\alpha(s) ds dt
\end{align*}
for any $x,y\in (0,\infty)$.

Therefore, by \cite[Lemma~3.6]{BDQS3}, there exists $C>0$ such that
\[\sup_{x\in I } \int_{(0,\infty)\setminus 2I} |\partial_x K_\phi^\alpha(x,y)|d\gamma_\alpha(y)\leq C\]
for every $I \in \mathcal{B}_1$, as claimed.

\bibliographystyle{acm}

\end{document}